\documentclass[10.5pt]{amsart}
\usepackage{amsmath,amsfonts,amssymb,amscd}
\usepackage[all]{xy}
\usepackage{graphics}
\usepackage{epsfig}
\usepackage{psfrag}

\def\ignore #1 {}

\newtheorem{lem}{Lemma}
\newtheorem{thm}{Theorem}

\newtheorem{dfn}{Definition}
\newtheorem{cor}{Corollary}
\newtheorem{rmk}{Remark}

\newtheorem{ex}{Example}

\def\del{\partial}

\begin{document}

\title{Gerbes and the Holomorphic Brauer Group of Complex Tori}
\author{Oren Ben-Bassat}
 \maketitle

\centerline{\bf }


{\tiny \tableofcontents }




\begin{abstract}
The purpose of this paper is to develop the theory of holomorphic gerbes on complex tori in a manner analogous to the classical theory for line bundles.  In contrast to past studies on this subject, we do not restrict to the case where these gerbes are torsion or  topologically trivial.  We give an Appell-Humbert type description of all
holomorphic gerbes on complex tori.  This gives an explicit, simple, cocycle
representative (and hence gerbe) for each equivalence class of holomorphic gerbes.
We also prove that a gerbe on the fiber product of four spaces over a common base is trivial as long as it is trivial upon restriction to any three out of the four spaces.
A fine moduli stack for gerbes on complex tori is constructed.  This involves the construction of a 'Poincar\'e' gerbe which plays a role analogous to the role of the Poincar\'e bundle in the case of line bundles.
\end{abstract}
\section{Introduction}\label{intro}  Giraud \cite{Gi1971} in 1971
began the study of {\it gerbes}, certain locally trivial stacks over a variety or scheme.   Equivalence
classes of gerbes are in one to one correspondence with the elements of a certain second cohomology group on the space.  One can define a gerbe over a space (or stack) as a torsor for $\mathcal{P}ic$, the stack of $\mathcal{O}^{\times}$ torsors on that space (or stack).  Brylinski includes some developments relating to gerbes in the analytic context in his book \cite{Br1993}.  Gerbes typically play the role of parameterizing certain moduli of geometric objects (\cite{Be2006}, \cite{BeBlPa2007}, \cite{DoGa2002}, \cite{DoPa2008}).  For a review on some things known about this group in the analytic context see \cite{Sc2005} and the references therein.  For
our purposes the crucial fact is that $\mathcal{O}^{\times}$ gerbes (with trivial band) on a complex manifold (or complex analytic space) $X$ are classified up to equivalence by
$H^{2}(X, \mathcal{O}^{\times})$.  By comparison, note that $\mathcal{O}^{\times}$ torsors (or their associated line bundles) are classified up to equivalence by $H^{1}(X, \mathcal{O}^{\times})$.  Previous studies (\cite{ElNa1983}, \cite{Be1972}, and \cite{Ho1972}) of the group $H^{2}(X, \mathcal{O}^{\times})$ for $X$ a complex torus or abelian variety have focused mainly on the representability by Azumaya algebras in the torsion case.  Line bundles on Abelian varieties (or complex tori) and their associated theta functions have a rich and distinguished history which we do not attempt to summarize here.  We begin in this paper a development for gerbes analogous to well known results for line bundles which can be found in a text-book on Abelian varieties or complex tori such as \cite{Po2003}, \cite{BiLa1999}, \cite{BiLa2004} or \cite{Mu1970}.  Based on results of Appell \cite{Ap1891} and
Humbert \cite{Hu1893} from the
1890s, Weil \cite{We1958} in 1958 gave the modern description of all line bundles on
complex tori known as the Appell-Humbert theorem.   We will prove a similar theorem giving in a concrete way a unique gerbe amongst the equivalence class corresponding to each element of $H^{2}(X, \mathcal{O}^{\times})$.  In contrast to the case of line bundles one cannot give a unique cocycle in this manner.  We will prove that in the relative setting that the assignment of a space $X$ to the group $H^{2}(X,\mathcal{O}^{\times})$ is a cubic functor in the sense of page 55 of \cite{Mu1970}.   We call this the theorem of the hyper-cube as it relates to gerbes on the product of four spaces.  In particular, this extends some results of Hoobler \cite{Ho1972} to the analytic setting.  We use this to derive formulas for the pullback of gerbes under various translation and multiplication maps, as well as isogenies.  We define the moduli stack for topologicaly trivial gerbes as $[H^{2}(X,\mathcal{O})/H^{2}(X,\mathbb{Z})]$ and define the universal (Poincar\'e) gerbe on the product of a torus and the moduli stack, proving that the moduli stack is fine.  In the second appendix and comments on future work, we call on some tools from and explain some relations with a recent preprint of Polishchuk \cite{Po2008} and comment on similarities with work of Felder, Henriques, Rossi, and Zhu \cite{FeHeRoZh2008}.
\section{Conventions}\label{convent}
In this paper $X$ will denote a complex torus of (complex) dimension $g$ written as $V/\Lambda$ where $\Lambda$ is a free group of rank $n=2g$ with $\Lambda \otimes \mathbb{R} = V$, a vector space with complex structure.  The quotient map will be denoted $p:V \to X$.  \\ The exponential map will be written without the $2\pi i$.  In other words
\[ \exp: \mathbb{C} \to \mathbb{C}^{\times} \]
always denotes the map
\[
 z \mapsto e^{2 \pi i z}.
 \]  Many of the computations in this paper make more explicit certain maps in the long exact sequence of cohomology groups coming from the following short exact sequence of sheaves of groups on $X$ induced by the exponential map:
 \begin{equation}\label{exponentialSES}
 0 \to \mathbb{Z} \to \mathcal{O}  \to \mathcal{O}^{\times} \to 1.
 \end{equation}
Unless stated explicitly, all vector spaces are over $\mathbb{R}$ or $\mathbb{C}$.  For any groups $\Gamma$ and $H$, the notation $\text{Alt}^{p}(\Gamma,H)$ denotes the maps from $\Gamma^{p}$ to $H$ that are skew-symmetric and group homomorphisms in each variable.  Here $\Gamma^{p}$ is the cartesian product of $p$ copies of $\Gamma$.  Although much of this paper concerns a complex torus, in sections \ref{mult} and \ref{univ} we will use {\em complex analytic spaces} \cite{De}, these are always assumed to be separated and reduced.  A reader who does not want to deal with this can feel free to substitute complex manifolds for complex analytic spaces.  In the relative context of a map $X \to S$ the same reader could assume that all fibers are complex manifolds.  For any space $M$ we denote by $\mathcal{P}ic(M)$ the groupoid of line bundles on $M$ and we denote by $\mathfrak{G}erbes(M)$ the $2-$category of gerbes on $M$.  This category has a monoidal structure which we denote by $\otimes$.  In particular we have the maps
 \[\mathcal{P}ic(M) \to \text{Pic}(M) = H^{1}(M, \mathcal{O}^{\times})
 \]
 denoted as
 \[\mathcal{L} \mapsto [\mathcal{L}]
 \]
 and
  \[\mathfrak{G}erbes(M) \to H^{2}(M, \mathcal{O}^{\times})
 \]
 denoted as
 \[\mathfrak{G} \mapsto [\mathfrak{G}].
 \]
 The reader is assumed to be familiar with group cohomology, we use the conventions of \cite{Mu1970}.   The boundary map in the cochain complex defining the group cohomology of a group $\Gamma$ acting on a module $R$ is denoted
 \begin{equation}\label{BoundaryMapGroupCoho}\delta: C^{p}(\Gamma,R) \to C^{p+1}(\Gamma,R).
 \end{equation}
\section{Acknowledgments}
Thank you to A. Polishchuk, P. Deligne, D. Kazhdan, T. Pantev, J. Block, J. de Jong, D. Gaitsgory, R. Hoobler, R. Livne, H. Farkas, M. Leyenson and many others for helpful and interesting conversations.  Thanks to anonymous referee for a helpful analysis.  Thanks to the mathematics department of the University of Salamanca and the Hebrew University where this work took place.

 \section{Cohomology and Alternating Classes}\label{coho}

Let the discrete group $\Gamma$ act freely and discontinuously on a complex analytic space $W$ and let $\mathcal{S}$ be a sheaf of groups on the quotient space $W/\Gamma$.  Let
 \[\rho:W \to W/\Gamma
 \]
be the quotient map.  Then there is a spectral sequence
\begin{equation}\label{eqn:SpectralSequence}
H^{p}(\Gamma, H^{q}(W,\rho^{-1}\mathcal{S})) \Longrightarrow H^{p+q}(W/\Gamma,\mathcal{S}).
\end{equation}
In particular, we get maps $H^{p}(\Gamma, H^{0}(W, \rho^{-1}\mathcal{S})) \to H^{p}(W/\Gamma, \mathcal{S})$ for all $p$.
In the case $p=1$, the resulting map $H^{1}(\Gamma, H^{0}(W, \rho^{-1}\mathcal{S})) \to H^{1}(W/\Gamma, \mathcal{S})$ is induced from the map sending a cocycle $\phi: \Gamma \to H^{0}(W,\rho^{-1}\mathcal{S})$ to the $\mathcal{S}$ torsor defined inside $\rho_{*} \mathcal{S}$ by the sections satisfying the equation $\gamma \cdot  \sigma = \phi(\gamma) \sigma$.  For any $p$, the map can be described in terms of a double complex comparing \v{C}ech and group cohomology, see \cite{Mu1970}.  The spectral sequence (\ref{eqn:SpectralSequence}) gives us maps
\begin{equation}\label{eqn:SpectralSequenceBoundary}
H^{p}(\Gamma, H^{0}(W,\rho^{-1}\mathcal{S})) \to H^{p}(W/\Gamma,\mathcal{S}).
\end{equation}
If we are also given sheaves of groups $\mathcal{T}$ and $\mathcal{R}$ and maps $\mathcal{S} \times \mathcal{T} \to \mathcal{R}$ which respect the group structure in each variable, the map (\ref{eqn:SpectralSequenceBoundary}) is compatible with the resulting cup products.
When the higher cohomology groups vanish, $H^{q}(W,\rho^{-1}\mathcal{S}) = 0$ for $q>0$, then the maps (\ref{eqn:SpectralSequenceBoundary}) become isomorphisms.

\begin{rmk}\label{rem:GerbeOnStackFromGroupCoho}
If we have any group action of $\Gamma$ on $W$, similar results hold as long as we pass to the stack $[W/\Gamma]$ as proven in Theorem A.6 of \cite{FeHeRoZh2008}.  In particular, if $H^{i}(W,\rho^{-1}\mathcal{S})) = 0$ for all $i>0$ then
\[H^{i}([W/\Gamma], \mathcal{S}) \cong H^{i}(\Gamma, H^{0}(W, \rho^{-1}\mathcal{S}))
\]
\end{rmk}
For the remainder of this section, we assume that $W$ is connected and contractible and $\Gamma$ is a free abelian group.
Let $R$ be a ring.  We need a skew-symmetrization map which we call $s$.

\begin{dfn}\label{SkewSymMap}
The skew-symmetrization map
\[s: \text{Map}(\Gamma^{p}, R) \to
\text{Map}(\Gamma^{p}, R)
\]
is defined by
\[s(f)(\gamma_{1}, \dots, \gamma_{p}) = \sum_{\sigma \in S_{p}}
(-1)^{\sigma}f(\gamma_{\sigma(1)}, \dots, \gamma_{\sigma(p)}).
\]
\end{dfn}
Note that when the ring $R$ is divisible the map $s$ has a section
\[ \frac{1}{p!}:Skew(\Gamma^{p},R) \to Map(\Gamma^{p}, R). \]
Our main case of interest will be $R=\mathbb{R}$ or $R=\mathbb{Z}$.
Let $\Gamma$ act on
$W$ to give a quotient $[W/\Gamma]$.  We let $\Gamma$ act trivially on the ring $R$, thought of as a constant sheaf on $W$.
\begin{lem}\label{lem:skew}
Consider the group cohomology $H^{p}(\Gamma,R)$.
The skew symmetrization map $s$ takes the group cocycles
$Z^{p}(\Gamma, R)$ to $\textup{Alt}^{p}(\Gamma,
R)$.  Furthermore, $s$ kills $B^{p}(\Gamma, R)$
and the resulting map from $H^{p}(\Gamma, R)$ to
$\textup{Alt}^{p}(\Gamma, R)$ is an isomorphism.
When composed with
\[\textup{Alt}^{p}(\Gamma, R)
\cong \wedge^{p}(\text{Hom}( \Gamma, R)) \cong
\wedge^{p}(H^{1}([W/\Gamma],R))\cong H^{p}([W/\Gamma],R)
\] it
agrees with the canonical isomorphism
\[H^{p}(\Gamma, R) \to H^{p}([W/\Gamma],R)
\]
which comes from the spectral sequence when $R$ is thought of as a
constant sheaf of groups on $[W/\Gamma]$.
\end{lem}

{\bf Proof.}

First of all suppose that $f \in B^{p}(\Gamma, R)$. This
means that we have a $g \in C^{p-1}(\Gamma,R)$  such
that
\begin{equation}
\begin{split}
f(\gamma_{0}, \dots, \gamma_{p-1}) & =  g(\gamma_{1}, \dots, \gamma_{p-1})
+\Big{[}\sum_{i=0}^{p-1}(-1)^{i+1}g(\gamma_{0}, \dots, \gamma_{i}
+\gamma_{i+1}, \cdots, \gamma_{p-1})\Big{]}  \\
& \quad + (-1)^{p}g(\gamma_{0},
\dots, \gamma_{p-2})
\end{split}
\end{equation}

When we skew-symmetrize $f$, one again has three terms: the
skew-symmetrizations of the first term, the term in the square
brackets, and the last term respectively.  Explicitly, we have
\begin{equation}
\begin{split}
 &\sum_{\sigma \in S_{p}} (-1)^{\sigma} g(\gamma_{\sigma (1)}, \dots, \gamma_{\sigma(p-1)}) \\
 & +\Big{[}\sum_{\sigma \in S_{p}} (-1)^{\sigma}\sum_{i=0}^{p-1}(-1)^{i+1}g(\gamma_{\sigma(0)}, \dots, \gamma_{\sigma(i)}
+\gamma_{\sigma(i+1)}, \cdots, \gamma_{\sigma(p-1)})\Big{]} \\
 & + \sum_{\sigma \in S_{p}}(-1)^{\sigma}(-1)^{p}g(\gamma_{\sigma(0)},
\dots, \gamma_{\sigma(p-2)}).
\end{split}
\end{equation}

The first and last of these new terms cancel because in the last
one, we can replace every appearance of $\sigma$ in the summand
with $\sigma \sigma'$ where $\sigma'$ is the permutation sending
$i$ to $i-1$ for $i \geq 0$ and $0$ to $p$. This permutation has
sign $(-1)^{(p-1)}$ making the first and last terms cancel.  For
the remaining terms we can replace $\sigma$ in the summands with
$\sigma\sigma'$ where $\sigma'$ is the permutation that flips $i$
and $i+1$.  This has sign $-1$ meaning the middle term is equal to
its own negative, and hence is zero.

Up to this point we have only used the fact that $R$ is an abelian
group.  Now we use the ring structure.
Notice that the
projection
\[Z^{p}(\Gamma, R) \to H^{p}(\Gamma, R) =
\cup^{p}H^{1}(\Gamma, R)
\]
has a section defined by
\begin{equation}\label{eqn:SkewSection}
[\alpha_{1}] \cup \cdots \cup [\alpha_{p}] \mapsto [(\gamma_{1},
\dots, \gamma_{p}) \mapsto \alpha_{1}(\gamma_{1}) \cdots
\alpha_{p}(\gamma_{p})].
\end{equation}
Therefore to prove that the
skew-symmetrization of a cocycle is multi-linear, it suffices to
observe that the skew-symmetrization of the image of $[\alpha_{1}]
\cup \cdots \cup [\alpha_{p}]$ under the section (\ref{eqn:SkewSection}) is
multi-linear.  This is clear.

To see that the map we have described agrees with the
map coming from the spectral sequence, it is enough to observe that both maps take cup products to
wedge products, and that the two maps agree for $p=0$ and
$p=1$.

\ \hfill $\Box$

In this paper, there are only two main cases where we will use the preceding part of this section.   The first case is $W=V$ and $\Gamma=\Lambda$, so that the stack $[W/\Gamma]$ is just a complex torus.  The second is the case $W = V \times \wedge^{2}\overline{V}^{\vee}$, $\Gamma = \Lambda \times \text{Alt}^{2}(\Lambda, \mathbb{Z})$.  The second case will be used in section \ref{univ}.   In the remainder of this section we will focus only on the first case.
The maps we have described are functorial with respect to maps of
the ring.  In particular, the relevant case will be the inclusion
$\mathbb{Z} \subset \mathbb{R}$.  In that case the situation is
summarized in the commutative diagram

\[
\xymatrix{\ar @{} [dr] |{}
  H^{n}(\Lambda,\mathbb{Z}) \ar[d]^s \ar[r] &  H^{n}(\Lambda,\mathbb{R})  \ar[r]^{=} & H^{n}(\Lambda,\mathbb{R}) \ar[d]^s  \\
  \text{Alt}^{n}(\Lambda, \mathbb{Z}) \ar[r]^{=} & \text{Alt}^{n}(\Lambda, \mathbb{Z}) \ar[u]^{\frac{1}{n!}} \ar[r] & \text{Alt}^{n}( \Lambda, \mathbb{R}) }
\]
where the outer vertical maps $s$ are isomorphisms.  If one picks a basis $\{\lambda^{i} \}$ of $\Lambda$ then for $n=2$ there is an explicit inverse $\sigma$ of the left most map $s$.  In fact, there is a group homomorphism
 \[\sigma: \textup{Alt}^{2}(\Lambda, \mathbb{Z}) \to Z^{2}(\Lambda, \mathbb{Z})
 \]
inducing a section to the map $s$.  It is defined by
\begin{equation}\label{eqn:DefOfSigma}\sigma(\mu)(\lambda_{1}, \lambda_{2})= \sum_{i<j}\mu(n_{1,i} \lambda^{i}, n_{2,j} \lambda^{j})
\end{equation}
where
\[\lambda_{\alpha} = \sum_{i} n_{\alpha,i} \lambda^{i}.
\]
\section{Recollection of Some Facts About Line Bundles}\label{LineAH}
In this section we recall the Appell-Humbert theorem which gives an explicit way to pick a line bundle out of each isomorphism class of line bundles on a complex torus.  We also review the Poincar\'{e} bundle and some of its properties.  Let $X=V/\Lambda$ be a complex torus.

Define $A(\Lambda)$ by
\[A(\Lambda) = \{E \in \text{Alt}^{2}(\Lambda,\mathbb{Z}) | E(ix,iy) = E(x,y) \} \subset \text{Alt}^{2}(\Lambda, \mathbb{Z}).
\]
Under the skew-symmetrization map $s$ preceded by the map coming from the spectral sequence
\begin{equation}
\label{eqn:blahequation}
H^{2}(X,\mathbb{Z}) \cong H^{2}(\Lambda, \mathbb{Z}) \cong \text{Alt}^{2}(\Lambda, \mathbb{Z}),
\end{equation}
$A(\Lambda)$ corresponds to the image of $H^{1}(X,\mathcal{O}^{\times})$ under the map coming from (\ref{exponentialSES}).  Notice that we have used the unique extension of $E$ to $\text{Alt}^{2}(V,\mathbb{R})$ where we denote it by the same letter.  We will do this kind of thing without explanation in the remainder of this article.  Furthermore notice that there is a natural identification
\begin{equation}\label{eqn:antiholOneform}
H^{1}(X, \mathcal{O}) \cong \overline{V}^{\vee}.
\end{equation}
Indeed, we have a Hodge projection map $H^{1}(X,\mathbb{R}) \to H^{1}(X,\mathcal{O})$ as well as
\begin{equation}\label{eqn:LineHodgeProj} H^{1}(X,\mathbb{R}) = \text{Hom}(\Lambda, \mathbb{R}) \subset V^{\vee} \oplus \overline{V}^{\vee} \twoheadrightarrow \overline{V}^{\vee}.
\end{equation}

Both these maps are isomorphisms so this gives us the isomorphism in (\ref{eqn:antiholOneform}).  Explicitly, the Hodge projection
\[\text{Hom}(\Lambda, \mathbb{C}) =H^{1}(\Lambda,\mathbb{C}) \to H^{1}(\Lambda, \mathcal{O}(V))
\]
is just
\[c \mapsto \Big{[}\lambda \mapsto \frac{c(\lambda)+ic(i\lambda)}{2} \Big{]}
\]
and the section of this map sends $l \in \overline{V}^{\vee}$ to the cohomology class coming from the element of $Z^{1}(\Lambda, \mathcal{O}(V))$ given by $[\lambda \mapsto 2Re(l(\lambda))]$.  We {\it define} $\Lambda^{\vee}$ as the image of $H^{1}(X, \mathbb{Z})$ in $\overline{V}^{\vee}$.  This comes out to be
\[\Lambda^{\vee} = \{l \in \overline{V}^{\vee} | 2 Re(l(\Lambda)) \subset \mathbb{Z}\}. \]
This is slightly different that the normal definition which substitutes $2 Re(l(\Lambda)) \subset \mathbb{Z}$ with the condition $Im(l(\Lambda)) \subset \mathbb{Z}$.  This discrepancy is discussed on page 87 of \cite{Mu1970}.  The two 'dual' tori that one gets from the two different possible definitions of $\Lambda^{\vee}$ are isomorphic.

In general we will denote the Hodge projection map induced by the inclusion $\mathbb{C} \to \mathcal{O}$ and given by either of the horizontal aarows in the following diagram
\[
\xymatrix{\ar @{} H^{p}(X, \mathbb{C}) \ar[r] \ar[d]^{\cong}& H^{p}(X, \mathcal{O})  \ar[d]^{\cong}\\
\text{Alt}^{p}(\Lambda, \mathbb{C}) \ar[r] & \wedge^{p} \overline{V}^{\vee}
}
\]
by
\begin{equation}\label{HodgeProjMap}
\Omega \mapsto \Omega^{H}.
\end{equation}
Notice that $H^{1}(X,\mathbb{Z})$ is isomorphic to $\Lambda^{\vee}$ via the map
\[Hom(\Lambda, \mathbb{Z})  \to \Lambda^{\vee}
\]
\[h \mapsto h^{H}
\]
which is an isomorphism compatible with equation (\ref{eqn:LineHodgeProj}) with inverse $\xi \mapsto \xi+\overline{\xi}$.
The following theorem (see e.g. \cite{Mu1970}) is implicit in the Appell-Humbert theorem as usually presented, however we chose to present it in this way to emphasize the analogy with our theorem in the case of gerbes.  In the following theorem we consider $(\overline{V}^{\vee}/\Lambda^{\vee}) \times A(\Lambda)$ as a category in the trivial way: the objects are points and the only morphisms are the identity maps.

\begin{thm}\label{lineAH}
There is a functor
\[\bold{ah}: (\overline{V}^{\vee}/\Lambda^{\vee}) \times A(\Lambda) \to \mathcal{P}ic(X)
\]
\[\bold{ah}([l], E) = \mathcal{L}_{([l], E)}
\]
such that $c_{1}(\mathcal{L}_{([l], E)})$ corresponds to $E$ under (\ref{eqn:blahequation}) and \[[\mathcal{L}_{([l], 0)}] = [l] \in H^{1}(X, \mathcal{O})/H^{1}(X,\mathbb{Z}) \subset H^{1}(X, \mathcal{O}^{\times})\]  using (\ref{eqn:antiholOneform}).  In the resulting isomorphism
\[(\overline{V}^{\vee}/\Lambda^{\vee}) \times A(\Lambda) \cong \textup{Pic}(X)
\]
the group structure on $(\overline{V}^{\vee}/\Lambda^{\vee}) \times A(\Lambda)$ induced from that on $\overline{V}^{\vee}$ and $\textup{Alt}^{2}(\Lambda, \mathbb{Z})$ corresponds to the tensor product of line bundles.
\end{thm}
{\bf Proof.}
Consider the homomorphism
\[A(\Lambda) \to Z^{1}(\Lambda, \mathcal{O}^{\times}(V))
\]
given by
\[E \mapsto \phi^{E}
\]
where
\[\phi^{E}_{\lambda}(v) =  \exp \Big{(}\vartheta(\lambda)
 -\frac{i}{2}E(iv,\lambda)+
 \frac{1}{2}E(v,\lambda)
 -\frac{i}{4}E(i\lambda,\lambda)\Big{)}.
\]

Here, the map \[\vartheta:\Lambda \to \mathbb{R},
\] sometimes called a {\it semi-character} for $E$,
is given in terms of some basis $\{ \lambda^{i} \}$ of $\Lambda$ as
\[\vartheta(\lambda) = \frac{1}{2} \sum_{i <j}E(n_{i} \lambda^{i}, n_{j} \lambda^{j})
\]
where $\lambda = \sum_{i} n_{i} \lambda^{i}$.

If we denote by $\mathcal{L}_{E}$ the line bundle corresponding in the sense of (\ref{eqn:Ob1CyclesToGroupoid})
to $\phi^{E}$ then its easy to check that $c_{1}(\mathcal{L}_{E}) = \del \phi^{E}$
corresponds under (\ref{eqn:blahequation}) to $E$.
Furthermore we can define a functor
\[F: [\overline{V}^{\vee} / \Lambda^{\vee}] \to  \mathcal{P}ic(X).
\]
On the level of objects, it sends $l$ to the line bundle $\mathcal{L}_{l}$
 corresponding in the sense of (\ref{eqn:Ob1CyclesToGroupoid}) to the constant cocycle in $Z^{1}(\Lambda, \mathcal{O}^{\times}(V))$ given by
\[\phi^{l}_{\lambda}(v) = \exp(l(\lambda)).
\]
The equivalence class $[\mathcal{L}_{l}] \in H^{1}(X,\mathcal{O}^{\times})$ comes from $l \in H^{1}(X,\mathcal{O})$ via the exponential map.
On the level of morphisms, suppose that $l_{2} - l_{1} = \xi \in \Lambda^{\vee}$.  Then we send $\xi$ to the isomorphism $\mathcal{L}_{l_{1}} \to \mathcal{L}_{l_{2}}$ determined in the sense of (\ref{eqn:Mor1CyclesToGroupoid}) by the element of $H^{0}(V,\mathcal{O}^{\times}) = C^{0}(\Lambda, \mathcal{O}^{\times}(V))$
\[v \mapsto \exp(-\overline{\xi(v)})
\]
which has boundary
\[\exp(-\overline{\xi(\lambda)})= \exp(\xi(\lambda)).
\]
$F$ determines a map \[
\underline{F}: (\overline{V}^{\vee} / \Lambda^{\vee}) \to \mathcal{P}ic(X),\] each element in $(\overline{V}^{\vee} / \Lambda^{\vee})$ is sent to the disjoint union of all images of $F$ of the primage of the element in  $\overline{V}^{\vee}$ modulo the obvious action of $\Lambda^{\vee}$.
Finally let
\[\bold{ah}(\alpha, E) = \underline{F}(\alpha) \otimes \mathcal{L}_{E}.
\]
\ \hfill $\Box$

There exists a Poincar\'{e} line bundle $\mathcal{P}$ on $X \times (\overline{V}^{\vee}/\Lambda^{\vee})$ with the property that under (\ref{eqn:antiholOneform}).
\[[\mathcal{P}|_{X \times \{\alpha\}}] = \alpha \in H^{1}(X,\mathcal{O})/H^{1}(X,\mathbb{Z})
\]
where addition in $\overline{V}^{\vee}$ corresponds to tensor product of line bundles.
We define the Poincar\'{e} bundle as the bundle corresponding (see section \ref{coho} and equation (\ref{eqn:Ob1CyclesToGroupoid})) to the cocycle
\[\psi \in Z^{1}(\Lambda, \mathcal{O}^{\times}(V))\]
given by
\[\psi_{\lambda, \xi}(v,l) = \exp\Big{(}l(\lambda) +\xi(\lambda) - \overline{\xi(v)}\Big{)}.
\]
One can replace the map $\underline{F}$ with the map
\[\alpha \mapsto \mathcal{P}|_{X \times \{\alpha\}}.
\]

The Poincar\'{e} bundle is needed to prove that $X^{\vee} = \overline{V}^{\vee}/\Lambda^{\vee}$ is a fine moduli space for topologically trivial line bundles on $X$ as formulated in the following theorem \cite{BiLa2004}.   In the cited book, they give a proof which works both in the algebraic and analytic context.  We give here an analytic proof.
\begin{thm}\label{LineFine} For any connected normal complex analytic space $T$ and any line bundle $\mathcal{L}$ on $X \times T$ such that $\mathcal{L}|_{X \times t}$ is topologically trivial for each $t$ there is a unique map $f:T \to X^{\vee}$ such that $(1,f)^{*} \mathcal{P} \cong \mathcal{L} \otimes \mathcal{C}$ where $\mathcal{C}$ is a line bundle trivial on each fiber $X \times \{t \}$.
\end{thm}
{\bf Proof.}
In fact the normality assumption is not needed.  Let $U_{T}$ be the universal cover of $T$.  Let $\rho:X \times T \to T$  and $\tilde{\rho}:X \times U_{T} \to U_{T}$ be the projection maps.   Consider the exact sequence
\[0 \to R^{1}\rho_{*} \mathbb{Z}\to R^{1} \rho_{*} \mathcal{O} \to (R^{1} \rho_{*} \mathcal{O}^{\times} )_{0} \to 1.
\]
The image of the class in $H^{0}(T, (R^{1} \rho_{*} \mathcal{O}^{\times} )_{0} )$ defined by restricting $\mathcal{L}$ to the fibers defines a class 
\[f_{\mathcal{L},\mathbb{Z}} \in \text{Hom}(\pi_{1}(T), \Lambda^{\vee}) \cong H^{1}(T, R^{1}\rho_{*} \mathbb{Z}).\]
Consider now the exact sequence
\[0 \to R^{1}\tilde{\rho}_{*} \mathbb{Z} \to R^{1} \tilde{\rho}_{*} \mathcal{O} \to (R^{1} \tilde{\rho}_{*} \mathcal{O}^{\times} )_{0} \to 1.
\]
Define \[f_{\mathcal{L},\mathbb{C}}:U_{T} \to \overline{V}^{\vee}\] to be any lift to $H^{0}(U_{T},R^{1}\tilde{\rho}_{*}\mathcal{O})$ of the class in $H^{0}(U_{T},(R^{1}\tilde{\rho}_{*}\mathcal{O}^{\times})_0)$ gotten by restricting the pullback of $\mathcal{L}$ to $X \times U_{T}$ to the fibers of the projection to $U_{T}$.  Such an element $f_{\mathcal{L},\mathbb{C}}$ exists because any obstructions live in
\[H^{1}(U_{T}, R^{1}\tilde{\rho}_{*} \mathbb{Z}) = H^{1}(U_{T}, Hom(\Lambda, \mathbb{Z}))=0.\]  The maps $f_{\mathcal{L},\mathbb{C}}$ and $f_{\mathcal{L},\mathbb{Z}}$ are compatible in the sense that for every $\chi \in \pi_{1}(T)$
\[\chi \cdot f_{\mathcal{L},\mathbb{C}} - f_{\mathcal{L},\mathbb{C}} = f_{\mathcal{L},\mathbb{Z}}(\chi).
\]
 Thus they together define a map \[f:T \cong U_{T}/\pi_{1}(T) \to \overline{V}^{\vee}/\Lambda^{\vee}= X^{\vee}.\]  The line bundle $(1,f)^{*}\mathcal{P}$ has a cocycle representative
\[(1,f)^{*}\psi \in
Z^{1}(\Lambda \times \pi_{1}(T),\mathcal{O}^{\times}(V \times U_{T})\] given by
\[(\lambda, \chi) \mapsto \psi_{\lambda,f_{\mathcal{L},\mathbb{Z}}(\chi)}(v,f_{\mathcal{L},\mathbb{C}}(u))
\]
where $\lambda \in \Lambda$, $\chi \in \pi_{1}(T)$, $v \in V$, and $u \in U_{T}$.  An apeal to the Leray spectral sequence which computes $H^{1}(X \times T, \mathcal{O}^{\times})$ via the projection to $T$ finishes the proof.

\ \hfill $\Box$

In the following two sections we will develop an analogous formalism in the case of gerbes, that is to say replacing $H^{1}(X,\mathcal{O}^{\times})$ with $H^{2}(X,\mathcal{O}^{\times})$.  In the next section (section \ref{TopClass}) we discuss the image in $H^{3}(X,\mathbb{Z}) \cong \text{Alt}^{3}(\Lambda, \mathbb{Z})$ and in the section after that (section \ref{GerbeAH}) we will lift these classes to cocycles and hence gerbes.  Here we prove our Appell-Humbert theorem for gerbes.  In section $\ref{univ}$ we find the analogue of the Poincar\'{e} sheaf, a universal gerbe parameterizing topologically trivial gerbes.  There we also prove we have constructed a fine moduli stack.

\section{Topological Classes of Gerbes}\label{TopClass}

Let $X=V/\Lambda$ be a complex torus.  From the short exact sequence
\[0 \to \mathbb{Z} \to \mathcal{O} \to \mathcal{O}^{\times} \to 1
\]
we know that the image of $H^{2}(X, \mathcal{O}^{\times})$ in
$H^{3}(X, \mathbb{Z})$ agrees with the kernel of the map of
$H^{3}(X, \mathbb{Z})$ to $H^{3}(X, \mathcal{O})$.  Since $X$ is a
K\"{a}hler manifold, the Hodge decomposition tells us that the image
of $H^{2}(X, \mathcal{O}^{\times})$ in $H^{3}(X, \mathbb{Z})$ is the kernel of the map
$H^{3}(X,\mathbb{Z}) \to H^{3}(X,\mathcal{O})$ which is
the following intersection in $H^{3}(X, \mathbb{C})$:
\[
H^{3}(X, \mathbb{Z}) \cap (H^{1,2}(X) \oplus H^{2,1}(X)).
\]
In this section we will give some information about the image of this map and for every element in the image.  The sub-variety of the moduli space of complex tori representing those tori which admit such a class is described by the equations
(\ref{htb_constraint2}) which we will derive.

We can give a parametrization of these elements as follows.
Consider the vector space $H^{3}(X,\mathbb{R})=
\wedge^{3}_{\mathbb{R}}\text{Hom}_{\mathbb{R}}(V,\mathbb{R})$.  Consider the inclusion map \[H^{1,2}(X) \oplus H^{2,1}(X) = (V^{\vee} \otimes \wedge^{2}
\overline{V}^{\vee}) \oplus (\wedge^{2} V^{\vee} \otimes
\overline{V}^{\vee}) \hookrightarrow \wedge^{3}(V \oplus \overline{V})^{\vee} =
H^{3}(X,\mathbb{C}).
\]
We will describe in an intrinsic way the real and integral elements of the image.  The real elements will be called $A(V)$ and the integral elements will be called $A(\Lambda)$.
\begin{dfn}
We define the holomorphic topological Brauer group $HTB(X)$ of a complex torus $X$ as
\[HTB(X) = \text{im}[H^{2}(X,\mathcal{O}^{\times}) \to H^{3}(X,\mathbb{Z})] = H^{3}(X,\mathbb{Z})
\cap (H^{1,2}(X) \oplus H^{2,1}(X)).\]
\end{dfn}
This is a free group and \[0 \leq rk(HTB(X)) \leq 2n\binom{n}{2}.\] because $dim_{\mathbb{R}} H^{1,2}(X) = 2n \binom{n}{2}$.
Consider the projection
\[p^{(1,2)+(2,1)}:  \wedge^{3}_{\mathbb{R}}\text{Hom}_{\mathbb{R}}(V,\mathbb{R}) \to
\text{Alt}^{3}(V,\mathbb{R})^{(1,2)+(2,1)}\] onto the $(1,2)+(2,1)$ part given by
\[(p^{(1,2)+(2,1)}(E))(x,y,z) =\frac{3}{4}E(x,y,z)+\frac{1}{4}( E(ix,iy,z)+E(x,iy,iz)+E(ix,y,iz)).
\]
\begin{dfn}\label{Adefinitions}
We define the subgroup $A(\Lambda)$ of
$\text{Alt}^{3}(\Lambda,\mathbb{Z})$
by
\[A(\Lambda) = \{ E \in
\text{Alt}^{3}(\Lambda,\mathbb{Z})
| p^{(1,2)+(2,1)}(E)=E
 \}
 \]
 and similarly
 \[A(V) = \{ E \in
\wedge^{3}_{\mathbb{R}}\text{Hom}_{\mathbb{R}}(V,\mathbb{R})
| p^{(1,2)+(2,1)}(E)=E
 \}.
 \]
\end{dfn}

For further use we describe the equation which defines both $A(\Lambda)$ and $A(V)$ as
\begin{equation}\label{str1}
E(x,y,z) = E(ix,iy,z)+E(x,iy,iz)+E(ix,y,iz)
\end{equation}
and its equivalent form
\begin{equation}\label{str2}
E(ix,iy,iz) = E(x,y,iz)+E(ix,y,z)+E(x,iy,z).
\end{equation}
Using the skew-symmetrization map we have the isomorphism
\[H^{3}(X,\mathbb{Z}) \cong H^{3}(\Lambda,\mathbb{Z}) \cong \text{Alt}^{3}(\Lambda, \mathbb{Z})
\]
which restricts to an isomorphism $HTB(X) \cong A(\Lambda)$.

As we will explain below, even if the Picard number $rk(NS(X))$ is zero, one could still have
 non-zero elements in $A(\Lambda)$.  This means that there exists
 tori with only topologically trivial line bundles, but which carry
 topologically non-trivial gerbes.

Recall that any complex torus $X$ is biholomorphic to the quotient
$V/\Lambda$ where we let $V= \mathbb{R}^{2g}$ with the complex structure
$J$ and $\Lambda = \Pi \mathbb{Z}^{2g}$.  Here $\Pi$ is an element of
$M(g \times 2g,\mathbb{C})$ thought of as a map
$\mathbb{R}^{2g} \to \mathbb{C}^{g}$ which satisfies
\begin{equation} \label{eqn:basic_thing}
i \Pi = \Pi \circ J
\end{equation}
as maps from $\mathbb{R}^{2g}$ to $\mathbb{C}^{g}$.
Let \[a_{i,j,k} = E(\Pi e_{i}, \Pi e_{j} , \Pi e_{k}) \] for the standard basis $e_{i}$ of $\mathbb{Z}^{2g}$.  Since $a$ is skew-symmetric, it is determined by its values $a_{i,j,k}$ where $i<j<k$.  The equation (\ref{str1}) reads
\[E(\Pi e_{i}, \Pi e_{j}, \Pi e_{k}) = E(i\Pi e_{i}, i \Pi e_{j}, \Pi e_{k}) +E(i\Pi e_{i},  \Pi e_{j}, i \Pi e_{k}) +E(\Pi e_{i}, i \Pi e_{j}, i\Pi e_{k}).
\]
Using (\ref{eqn:basic_thing}) this equation becomes
\begin{equation}\label{htb_constraint}
a_{i,j,k} = a_{l,m,k} J_{l,i}J_{m,j}+ a_{i,m,n}J_{m,j}J_{n,k}+a_{l,j,n}J_{l,i}J_{n,k}
\end{equation}
for $Je_{s} = J_{t,s}e_{t}$ and $a \in \mathbb{Z}^{\binom{2g}{3}}.$  We will reexpress condition as an intersection of principal divisors on the moduli space of complex tori, each divisor corresponding to a holomorphic function of the parameters $\tau$.

Consider the matrix
\[\zeta = \left( \begin{array}{c}
1 \\
-\tau \end{array} \right)\\
\]

Then $J \zeta = -i \zeta$.  So we have $J_{z,p} \zeta_{p,q} = -i \zeta_{z,q}$.  In order to rewrite  equation (\ref{htb_constraint}) in terms of complex structure parameters $\tau$, we multiply (\ref{htb_constraint}) by $\zeta_{i,s} \zeta_{j,t} \zeta_{k,u}$ and sum over repeated indices.  Notice that there is nothing lost by restricting to the case $s < t< u$.   After some simplifications to the right hand side we get
\[a_{i,j,k} \zeta_{i,s} \zeta_{j,t} \zeta_{k,u} = - 3 a_{i,j,k} \zeta_{i,s} \zeta_{j,t} \zeta_{k,u}.
\]
Therefore the condition (\ref{htb_constraint}) implies
\[a_{i,j,k} \zeta_{i,s} \zeta_{j,t} \zeta_{k,u} = 0.
\]
When we expand this out, we get
\begin{equation}\label{htb_constraint2}
\begin{split}
 & a_{s,t,u} -a_{p+g,t,u} \tau_{p,s} -a_{s,q+g,u} \tau_{q,t} -a_{s,t,r+g} \tau_{r,u} \\
 & + a_{s,q+g,r+g} \tau_{q,t} \tau_{r,u} +a_{p+g,t,r+g} \tau_{p,s} \tau_{r,u} +a_{p+g,q+g,u}\tau_{p,s} \tau_{q,t} \\
 & - a_{p+g,q+g,r+g} \tau_{p,s} \tau_{q,t} \tau_{r,u} = 0.
\end{split}
\end{equation}
The number of solutions to equation (\ref{htb_constraint2}) could possibly be maximized by choosing the real and imaginary parts of the entries
of $\tau$ to be rational numbers.  Equation (\ref{htb_constraint2}) has no solutions if the real and imaginary parts of the entries of $\tau$ are chosen to be algebraically independent over $\mathbb{Q}$.

In the special case of a complex torus of complex dimension $3$, equation (\ref{htb_constraint2})  becomes

\begin{equation}\label{dim3htb_constraint2}
\begin{split}
 & a_{0,1,2}-\sum_{p}a_{1,2,p+3} \tau_{p,0} + \sum_{q}a_{0,2,q+3} \tau_{q,1} -\sum_{r}a_{0,1,r+3} \tau_{r,2} \\
& + \sum_{q<r}a_{0,q+3,r+3} (\tau_{q,1} \tau_{r,2}-\tau_{r,1} \tau_{q,2})- \sum_{p<r}a_{1,p+3,r+3} (\tau_{p,0} \tau_{r,2}-\tau_{r,0} \tau_{p,2}) \\
& + \sum_{p<q}a_{2,p+3,q+3}(\tau_{p,0} \tau_{q,1}-\tau_{q,0} \tau_{p,1}) \\
& - a_{3,4,5} \det(\tau) = 0.
\end{split}
\end{equation}
\begin{lem}
Equation (\ref{htb_constraint2}) is actually equivalent to equation (\ref{str1}).
\end{lem}
{\bf Proof.}
Notice that $\left( \begin{array}{cc}
\zeta & \overline{\zeta} \end{array} \right)$ in an invertible matrix and so the equation (\ref{htb_constraint}) is equivalent to the {\it contracted equation} that we get by multiplying both sides of equation (\ref{htb_constraint}) by
\[\left( \begin{array}{cc}
\zeta & \overline{\zeta} \end{array} \right)_{i,s} \left( \begin{array}{cc}
\zeta & \overline{\zeta} \end{array} \right)_{j,t} \left( \begin{array}{cc}
\zeta & \overline{\zeta} \end{array} \right)_{k,u}
\]
and summing over $i$, $j$, and $k$.
Notice that $J\left( \begin{array}{cc}
\zeta & \overline{\zeta} \end{array} \right) = \left( \begin{array}{cc}
-i \zeta & i \overline{\zeta} \end{array} \right)$ and so most of the equations from the expansion come out to put no constraints on $a$.  Indeed, the terms with two copies of $\zeta$ and one copy of $\overline{\zeta}$ place no constraints on $a$.  To see this notice that for the terms with two copies of $\zeta$ and one copy of $\overline{\zeta}$ the right hand side of equation (\ref{htb_constraint}) contributes three terms to the contracted equation: two of them equal to the left hand side of the contracted equation, and one equal to minus the left hand side of the contracted equation.  The same thing happens with the terms with one copy of $\zeta$ and two copies of $\overline{\zeta}$. The only remaining equations are those with three copies of $\zeta$, or three copies of $\overline{\zeta}$.  But we have already accounted for the former in equation (\ref{htb_constraint2}), and the later comes about from the former by complex conjugation.
\ \hfill $\Box$
\begin{rmk}
For $g=1$ all holomorphic gerbes are trivial and for $g=2$ the equation (\ref{htb_constraint2}) is a trivial equation as explained further in example \ref{stuff}.  For $g=3$ (and for all $g>3$ as well) the equation might have no solutions as can be seen by choosing $\tau$ to be some matrix of complex numbers with $\det(Im(\tau)) \neq 0$ whose collection of real and imaginary parts are a set of real numbers algebraically independent over the rational numbers.  This is analogous to the fact that the generic torus of complex dimension greater than or equal to $2$ has only topologically trivial line bundles.
\end{rmk}

\begin{ex}\label{stuff}

Suppose that $X$ is a complex torus of dimension $g=2$.  The Neron Severi
group could be trivial but nevertheless, the holomorphic topological Brauer
group is never trivial.  Indeed in the case $g=2$ we have
$HTB(X) = H^{3}(X,\mathbb{Z}) \cong \mathbb{Z}^{4}$ since
$H^{3,0}(X) = H^{0,3}(X) = (0)$.  This is analogous to the fact that
every elliptic curve $C$ over $\mathbb{C}$ has a non-trivial
Neron-Severi group $H^{1,1}(C,\mathbb{Z})
= H^{2}(C,\mathbb{Z})$.  This tells us
that (\ref{htb_constraint}) is satisfied for any choice
of the coefficients $a$.  For instance, taking a basis
$\{ e_{0}, e_{1}, e_{2}, e_{3} \}$ and $a_{0,1,2} = 1$,
and all the other entries zero when possible this says that
any complex structure $J$ on a $4$ dimensional real vector space satisfies
\[J_{0,0}J_{1,1}+J_{1,1}J_{2,2}+J_{0,0}J_{2,2}
-J_{0,1}J_{1,0}-J_{1,2}J_{2,1}-J_{0,2}J_{2,0} = 1
\]
in any basis.

The analysis in this example proves that the map
\[H^{1}(X, \mathcal{O}^{\times}) \otimes H^{1}(X, \mathbb{Z}) \to H^{2}(X, \mathcal{O}^{\times})
\]
induced by the map
\[\mathcal{O}^{\times} \times \mathbb{Z} \to \mathcal{O}^{\times}\]
given by
\[(f,n) \mapsto f^{n}
\]
cannot be surjective in general.
\end{ex}

\begin{ex}
In this example we look at three dimensional complex tori $X$ which are the product of three elliptic curves, each of which has purely imaginary period.  We show that the group $HTB(X)$ can discriminate amongst different such tori even amongst tori with fixed Picard number.  Notice also that algebraic dimension of the product of three elliptic curves is always $3$, showing that the group $HTB(X)$ has more information than the invariants which are typically studied.
Let $X = \mathbb{R}^{6} / \Pi \mathbb{Z}^{6}$
where $\Pi = (\tau, 1_{3})$ as is always possible.  Here $\tau \in M_{3}(\mathbb{C})$
is a three by three matrix with the imaginary part non-degenerate and $1_{3}$
is the three by three identity matrix.  The complex structure is given by the formula
\[J = \left( \begin{array}{cc}
y^{-1}x & y^{-1} \\
-y-xy^{-1}x & -xy^{-1}  \end{array} \right)\\
\]
where
\[x = Re(\tau)
\]
and
\[y = Im(\tau).
\]
Since we are looking at the product of three elliptic curves with purely imaginary periods, we consider complex structures of the form
\[\tau = \left( \begin{array}{ccc}
i\alpha & 0 & 0 \\
0 & i\beta & 0 \\
0 & 0 & i\gamma  \end{array} \right)\\
\]
where $\alpha, \beta, \gamma \in \mathbb{R}-0$.

Because $\tau$ is diagonal, the equation (\ref{dim3htb_constraint2}) becomes
\begin{equation}\label{eqn:diag_reduc}
\begin{split}
& a_{0,1,2}-a_{1,2,3} \tau_{0,0}+a_{0,2,4}\tau_{1,1}-a_{0,1,5}\tau_{2,2} \\
& +a_{0,4,5}\tau_{1,1}\tau_{2,2} -a_{1,3,5}\tau_{0,0}\tau_{2,2}+a_{2,3,4}\tau_{0,0} \tau_{1,1} \\
& -a_{3,4,5} \tau_{0,0} \tau_{1,1} \tau_{2,2} = 0.
\end{split}
\end{equation}

Due to the fact that $\tau$ is purely imaginary, equation (\ref{eqn:diag_reduc}) breaks up into the following two equations
\begin{equation}
\begin{split}
& a_{0,1,2} -a_{0,4,5}\beta \gamma +a_{1,3,5}\alpha\gamma - a_{2,3,4} \alpha \beta  = 0.
\end{split}
\end{equation}
and
\begin{equation}
\begin{split}
-a_{1,2,3} \alpha+a_{0,2,4}\beta-a_{0,1,5}\gamma +a_{3,4,5} \alpha \beta \gamma = 0.
\end{split}
\end{equation}
Let $R$ be the rank of the solution space to the above two equations.
The $12$ integers
\[a_{0,1,3}, a_{0,1,4},  a_{0,2,3}, a_{0,2,5}, a_{0,3,4},a_{0,3,5}, a_{1,2,4}, a_{1,2,5}, a_{1,3,4}, a_{1,4,5}, a_{2,3,5} \ \ \ \ \text{and} \ \ \ \  a_{2,4,5}\] do not appear in these equations and are therefore unconstrained.
Thus
\[rk(HTB(X(\alpha, \beta, \gamma))) = 12 + R.
\]
An easy computation using Proposition 3.4 on page 10 of \cite{BiLa1999} shows that the Neron-Severi group has rank
\[rk(NS(X(\alpha, \beta, \gamma))) = 3 + R_{1} + R_{2}
\]
where
\[R_{1} = \# \Big{(}\{\frac{\alpha}{\beta}, \frac{\beta}{\gamma}, \frac{\alpha}{\gamma} \} \cap \mathbb{Q}\Big{)}
\]
and
\[R_{2} = \# \Big{(}\{ \alpha \beta, \alpha \gamma, \beta \gamma \} \cap \mathbb{Q}\Big{)}.
\]

Therefore
\[\alpha=1, \beta = \sqrt{2},\gamma = \sqrt{3} \Longrightarrow  rk(NS(X(\alpha, \beta, \gamma))) = 3, rk(HTB(X(\alpha, \beta, \gamma))) = 12\]

while in contrast
\[\alpha = 1, \beta = \sqrt{2}, \gamma = \frac{\beta}{1-\beta}  \Longrightarrow  rk(NS(X(\alpha, \beta, \gamma))) = 3, rk(HTB(X(\alpha, \beta, \gamma))) \geq 14.\]

Indeed, two new solutions for $a$ are
\[a_{0,4,5} = a_{1,3,5} = a_{2,3,4} = 1 \]
\[
a_{i,j,k} = 0 \ \ \ \ \text{for other values of}\ \ \ \ \{i,j,k\}
\]
and
\[a_{0,2,4} = a_{0,1,5} = a_{3,4,5} = 1
\]
\[
a_{i,j,k} = 0 \ \ \ \ \text{for other values of} \ \ \ \ \{i,j,k\}.
\]
If we chose instead $\alpha \in \mathbb{Q}$, $\beta \in \mathbb{Q}$, and $\gamma \in \mathbb{Q}$ then both the rank of $NS(X)$ and the rank of $HTB(X)$ are maximal, being $9$ and $18$ respectively.


\end{ex}
We now explain more carefully the nature of the projections onto the $(1,2)+(2,1)$ part which we have used above.  Recall the definition of $A(V)$ from Definition \ref{str1}.  We have an isomorphism
\begin{equation}\label{AofV}\left( V^{\vee} \otimes \wedge^{2}
\overline{V}^{\vee} \bigoplus  \wedge^{2}
V^{\vee} \otimes \overline{V}^{\vee}\right)_{\mathbb{R}}  \cong A(V)
\end{equation}
induced by the canonical isomorphism
\[\wedge^{3}V^{\vee} \bigoplus V^{\vee} \otimes \wedge^{2}
\overline{V}^{\vee} \bigoplus   \wedge^{2} V^{\vee} \otimes
\overline{V}^{\vee} \bigoplus \wedge^{3}\overline{V}^{\vee}=
\wedge^{3}\left(V^{\vee} \bigoplus \overline{V}^{\vee}\right) =
\wedge^{3}_{\mathbb{C}}(\text{Hom}_{\mathbb{R}}(V,\mathbb{C}))
\]
Using the inclusions
\[\left(V^{\vee} \otimes \wedge^{2}
\overline{V}^{\vee} \bigoplus   \wedge^{2} V^{\vee} \otimes
\overline{V}^{\vee}\right)_{\mathbb{R}} \subset
\text{Hom}_{\mathbb{R}}(V \otimes_{\mathbb{R}} V
\otimes_{\mathbb{R}}V, \mathbb{C})
\] and
\[
\wedge^{3}_{\mathbb{R}}\text{Hom}_{\mathbb{R}}(V,\mathbb{R})
\subset \text{Hom}_{\mathbb{R}}(V \otimes_{\mathbb{R}} V
\otimes_{\mathbb{R}}V, \mathbb{C})
\]
Recall the definition of $s$ from equation (\ref{SkewSymMap}).  The maps giving the isomorphism in equation (\ref{AofV}) are
\[\frac{1}{3!}s: \left( V^{\vee} \otimes \wedge^{2}
\overline{V}^{\vee} \bigoplus   \wedge^{2} V^{\vee} \otimes
\overline{V}^{\vee}\right)_{\mathbb{R}}  \to A(V).
\]
and its inverse
\[A(V) \to \left(V^{\vee} \otimes \wedge^{2}
\overline{V}^{\vee} \bigoplus   \wedge^{2}
V^{\vee} \otimes \overline{V}^{\vee}\right)_{\mathbb{R}}
\]
is given by
\[E \mapsto \Bigg{[}(x,y,z) \mapsto \frac{3}{4}\Big{(}E(x,y,z)+E(ix,iy,z)
+E(ix,y,iz)-E(x,iy,iz)\Big{)}\Bigg{]}.\]

The expression appearing on the right hand side of the above equation is the real part of the following element of $V^{\vee} \otimes \wedge^{2} \overline{V}^{\vee}:$
\begin{eqnarray*}
& & \frac{3}{4}\Big{(}E(x,y,z)+E(ix,iy,z)+E(ix,y,iz)-E(x,iy,iz) \\
& & +i(E(x,y,iz)+E(x,iy,z)-E(ix,y,z)+E(ix,iy,iz))\Big{)}
\end{eqnarray*}
 The above element is the analogue in this gerbey context to the Hermitian form living in $V^{\vee} \otimes \overline{V}^{\vee}$ which corresponds to the first Chern class of a line bundle.  The real part of this Hermetian form corresponds to the $E$ mentioned for the line bundles case in Section \ref{LineAH}.  However, in our analysis of gerbes it is $E \in A(\Lambda)$ that will play the fundamental role.

\section{The Appell-Humbert Theorem for Gerbes}\label{GerbeAH}
In the following, we will construct elements of $Z^{2}(\Lambda,
\mathcal{O}^{\times}(V))$ as the image of elements of
$\text{Map}(\Lambda \times \Lambda , \mathcal{O}(V))$ under the
exponential map. Clearly this works if and only if the boundary of
the element is integral and so it will be killed by the exponential. In general
looking at the diagram
\[\begin{CD}
\text{Map}(\Lambda^{p+1}, \mathcal{O}(V)) \\
@AA{\delta}A \\
\text{Map}(\Lambda^{p}, \mathcal{O}(V)) @>\exp>>
\text{Map}(\Lambda^{p},\mathcal{O}^{\times}(V))
\end{CD}
\]

where the definition of $\delta$ can be found in (\ref{BoundaryMapGroupCoho}).  We can clearly see that

\[\delta^{-1}(\text{Map}(\Lambda^{p+1}, \mathbb{Z})) =
\exp^{-1}\Big{(}Z^{p}(\Lambda,\mathcal{O}^{\times}(V))\Big{)}.
\]
We would like to
find a cocycle $\Phi^{E} = \exp(\Theta^{E}) \in Z^{2}(\Lambda, \mathcal{O}^{\times}(V))$
which satisfies $\delta \exp(\Theta^{E}) = E$.  We do this in two stages.
First, we describe the preimage of $E$ under the map $s \circ \del
\circ \exp$:
\[
\xymatrix{ \text{Map}(\Lambda^{2}, \mathcal{O}(V)) \ar[r]^{exp} & \text{Map}(\Lambda^{2}, \mathcal{O}^{\times}(V)) \ar[r]^{\del} & \text{Map}(\Lambda^{3}, \mathbb{Z})  \ar[r]^{s} & \text{Map}(\Lambda^{3}, \mathbb{Z}).
}\]
Of course the standard map $\del$ is a boundary map in the long exact sequence in group cohomology coming from the exponential short exact sequence. It is induced by a map also called $\del$ from $Z^{2}(\Lambda, \mathcal{O}^{\times}(V))$ to $Z^{3}(\Lambda, \mathbb{Z})$ which takes  $B^{2}(\Lambda, \mathcal{O}^{\times}(V))$ to $B^{3}(\Lambda, \mathbb{Z})$.   Here, however, we simply use the same formula and symbol to define a map  $\text{Map}(\Lambda^{2}, \mathcal{O}^{\times}(V))$ to $\text{Map}(\Lambda^{3}, \mathbb{Z})$.
Therefore, in the first step, we find solutions $\Theta^{E}$ to the equation
\begin{equation}\label{eqn:FirstStepEquation}
s \circ \del \circ \exp (\Theta^{E}) = s \circ \delta (\Theta^{E}) = E.
\end{equation}
 Second, we find an
element in this preimage such that when we apply $\delta$ the
result is integral. In
other words, we impose the further constraint that
\begin{equation}\label{eqn:SecondStepEquation}
\delta \Theta^{E} \in \text{Map}(\Lambda \times \Lambda \times \Lambda,
\mathbb{Z}).
\end{equation}  The exponential of this element, $\Phi^{E} = \exp(\Theta^{E})$, is then the
desired element of $Z^{2}(\Lambda, \mathcal{O}^{\times}(V))$.

Let $\Theta^{E} \in \text{Map}(\Lambda \times \Lambda, \mathcal{O}(V))$
be given by \[(\lambda_{1},\lambda_{2}) \mapsto \Theta^{E}_{\lambda_{1},\lambda_{2}}(v)\] where
\begin{equation}\label{eqn:HplusBeta}
\Theta^{E}_{\lambda_{1},\lambda_{2}}(v) = H_{\lambda_{1},\lambda_{2}}(v)
+\beta_{\lambda_{1},\lambda_{2}}
\end{equation}

here we take
\begin{equation}\label{eqn:Heqn}
\begin{split}
H_{\lambda_{1},\lambda_{2}}(v) &= \frac{1}{8}\Big{(}E(v,\lambda_{1},\lambda_{2})+\frac{1}{2}E(iv,i\lambda_{1},\lambda_{2})
+\frac{1}{2}E(iv,\lambda_{1},i\lambda_{2})\Big{)} \\
&  \quad +\frac{i}{8}\Big{(}\frac{1}{2}E(v,i\lambda_{1},\lambda_{2})+\frac{1}{2}E(v,\lambda_{1},i\lambda_{2})- E(iv,\lambda_{1},\lambda_{2})\Big{)}
\end{split}
\end{equation}
and $\beta_{\lambda_{1},\lambda_{2}}$ are (for now) arbitrary complex constants.  Notice that $\Theta^{E}_{\lambda_{1},\lambda_{2}}$ is holomorphic because \[H_{\lambda_{1},\lambda_{2}}(iv) = i H_{\lambda_{1},\lambda_{2}}(v).\]
Also we claim $\Theta^{E}$ is
a solution to equation (\ref{eqn:FirstStepEquation}).  Indeed
using the additivity of $E$ in its entries, one has
\begin{equation}
\begin{split}
(\delta H)_{\lambda_{1},\lambda_{2},\lambda_{3}}(v) & =  H_{\lambda_{2},\lambda_{3}}(v+\lambda_{1}) - H_{\lambda_{1}+\lambda_{2},\lambda_{3}}(v)+H_{\lambda_{1},\lambda_{2}+\lambda_{3}}(v)-H_{\lambda_{1},\lambda_{2}}(v) \\
& = H_{\lambda_{2},\lambda_{3}}(\lambda_{1}).
\end{split}
\end{equation}
Thus $s\circ \delta (\psi)$ is the skew symmetrization of the map
\begin{equation}
\begin{split}
(\lambda_{1},\lambda_{2},\lambda_{3}) \mapsto &  H_{\lambda_{2},\lambda_{3}}(\lambda_{1}) \\
&   =  \frac{1}{8}\Big{(}E(\lambda_{1},\lambda_{2},\lambda_{3})+\frac{1}{2}E(i\lambda_{1},i\lambda_{2},\lambda_{3})
+\frac{1}{2}E(i\lambda_{1},\lambda_{2},i\lambda_{3})\Big{)} \\
&  \quad + \frac{i}{8}\Big{(}\frac{1}{2}E(\lambda_{1},i\lambda_{2},\lambda_{3})+\frac{1}{2}E(\lambda_{1},\lambda_{2},i\lambda_{3})- E(i\lambda_{1},\lambda_{2},\lambda_{3})\Big{)}.
\end{split}
\end{equation}
Define
\begin{equation}\label{eqn:keqn}
k(\lambda_{1},\lambda_{2},\lambda_{3}) = \text{Re}(H_{\lambda_{2},\lambda_{3}}(\lambda_{1}))=  \frac{1}{8}\Big{(}E(\lambda_{1},\lambda_{2},\lambda_{3})+\frac{1}{2}E(i\lambda_{1},i\lambda_{2},\lambda_{3})
+\frac{1}{2}E(i\lambda_{1},\lambda_{2},i\lambda_{3})\Big{)}
\end{equation}
and
\begin{equation}\label{eqn:leqn}
l(\lambda_{1},\lambda_{2},\lambda_{3}) = \text{Im}(H_{\lambda_{2},\lambda_{3}}(\lambda_{1})) = \frac{1}{8}\Big{(}\frac{1}{2}E(\lambda_{1},i\lambda_{2},\lambda_{3})
+\frac{1}{2}E(\lambda_{1},\lambda_{2},i\lambda_{3})- E(i\lambda_{1},\lambda_{2},\lambda_{3})\Big{)}.
\end{equation}
The skew-symmetrization of the real part gives us
\begin{equation} \label{eqn:skewK}
\begin{split}
 & (sk)(\lambda_{1},\lambda_{2},\lambda_{3}) \\
 & = \frac{1}{8} \bigg{(} 6E(\lambda_{1},\lambda_{2},\lambda_{3})+\frac{1}{2}E(i\lambda_{1},i\lambda_{2},\lambda_{3})+
\frac{1}{2}E(i\lambda_{1},\lambda_{2},i\lambda_{3}) - \frac{1}{2}E(i\lambda_{1},i\lambda_{3},\lambda_{2})-\frac{1}{2}E(i\lambda_{1},\lambda_{3},i\lambda_{2}) \\
 & \quad -\frac{1}{2}E(i\lambda_{2},i\lambda_{1},\lambda_{3})-\frac{1}{2}E(i\lambda_{2},\lambda_{1},i\lambda_{3}) +\frac{1}{2}E(i\lambda_{2},i\lambda_{3},\lambda_{1})+\frac{1}{2}E(i\lambda_{2},\lambda_{3},i\lambda_{1})  \\
 & \quad +\frac{1}{2}E(i\lambda_{3},i\lambda_{1},\lambda_{2})+\frac{1}{2}E(i\lambda_{3},\lambda_{1},i\lambda_{2}) -\frac{1}{2}E(i\lambda_{3},i\lambda_{2},\lambda_{1}) - \frac{1}{2}E(i\lambda_{3},\lambda_{2},i\lambda_{1}) \bigg{)} \\
 & =  \frac{1}{8}\bigg{(}6E(\lambda_{1},\lambda_{2},\lambda_{3})
 +2\Big{(}E(i\lambda_{1},i\lambda_{2},\lambda_{3})+E(i\lambda_{1},\lambda_{2},i\lambda_{3})+E(\lambda_{1},i\lambda_{2},i\lambda_{3})\Big{)}\bigg{)} \\
 & = E(\lambda_{1},\lambda_{2},\lambda_{3})
\end{split}
\end{equation}
where at the last step we have used equation (\ref{str2}). The
skew-symmetrization of the imaginary part comes out to be
zero.  Indeed

\begin{equation}
\begin{split}
 & 8(sl)(\lambda_{1},\lambda_{2},\lambda_{3}) \\
 & = E(i\lambda_{1},\lambda_{2},\lambda_{3})
 -E(i\lambda_{1},\lambda_{3},\lambda_{2})
 -E(i\lambda_{2},\lambda_{1},\lambda_{3})
 +E(i\lambda_{2},\lambda_{3},\lambda_{1})
 +E(i\lambda_{3},\lambda_{1},\lambda_{2})
 -E(i\lambda_{3},\lambda_{2},\lambda_{1}) \\
 &
  -\frac{1}{2}\Big{(}E(\lambda_{1},i\lambda_{2},\lambda_{3})
 -E(\lambda_{1},i\lambda_{3},\lambda_{2})
 -E(\lambda_{2},i\lambda_{1},\lambda_{3})
 +E(\lambda_{2},i\lambda_{3},\lambda_{1})
 +E(\lambda_{3},i\lambda_{1},\lambda_{2})
 -E(\lambda_{3},i\lambda_{2},\lambda_{1})\Big{)} \\
&
-\frac{1}{2}\Big{(}E(\lambda_{1},\lambda_{2},i\lambda_{3})
-E(\lambda_{1},\lambda_{3},i\lambda_{2})
-E(\lambda_{2},\lambda_{1},i\lambda_{3})
+E(\lambda_{2},\lambda_{3},i\lambda_{1})
+E(\lambda_{3},\lambda_{1},i\lambda_{2})
-E(\lambda_{3},\lambda_{2},i\lambda_{1})\Big{)} \\
 & = 2E(i\lambda_{1},\lambda_{2},\lambda_{3})
 +2E(\lambda_{1},i\lambda_{2},\lambda_{3})
 +2E(\lambda_{1},\lambda_{2},i\lambda_{3}) \\
 &
 -\frac{1}{2}\Big{(}2E(\lambda_{1},i\lambda_{2},\lambda_{3})
 +2E(\lambda_{1},\lambda_{2},i\lambda_{3})
 +2E(i\lambda_{1},\lambda_{2},\lambda_{3})\Big{)} \\
 &
 -\frac{1}{2}\Big{(}2E(\lambda_{1},\lambda_{2},i\lambda_{3})
 +2E(\lambda_{1},i\lambda_{2},\lambda_{3})
 +2E(i\lambda_{1},\lambda_{2},\lambda_{3})\Big{)} \\
 & = 0
\end{split}
\end{equation}

Using Lemma \ref{lem:skew} have
\[s(\delta(\beta)) = 0.
\]
Recalling equation (\ref{eqn:HplusBeta}) we can conclude that equation (\ref{eqn:FirstStepEquation}) holds.

This concludes the first step.
For the second step we need to consider the integrality equation (\ref{eqn:SecondStepEquation}) which
reads
\begin{equation}\label{eqn:integral}
H_{\lambda_{2},\lambda_{3}}(\lambda_{1}) + \beta_{\lambda_{2},\lambda_{3}}-\beta_{\lambda_{1}
+\lambda_{2},\lambda_{3}}+\beta_{\lambda_{1},\lambda_{2}
+\lambda_{3}}-\beta_{\lambda_{1},\lambda_{2}}
 \in \mathbb{Z}.
\end{equation}
for all $\lambda_{1},\lambda_{2},\lambda_{3} \in \Lambda$.  Let
\[\beta = \beta'+ i \beta''
\]
be the decomposition of $\beta$ into real and imaginary parts.   Consider
the decomposition of equation (\ref{eqn:integral}) into real and
imaginary parts.  In order to satisfy the imaginary part we need
\begin{equation}\label{kill_imag_part}
(\delta \beta'')_{\lambda_{1},\lambda_{2},\lambda_{3}} +l(\lambda_{1},\lambda_{2},\lambda_{3}) = 0
\end{equation}
or
\[
\beta''_{\lambda_{2},\lambda_{3}}-\beta''_{\lambda_{1}+\lambda_{2},\lambda_{3}}+
\beta''_{\lambda_{1},\lambda_{2}+\lambda_{3}}
-\beta''_{\lambda_{1},\lambda_{2}}+l(\lambda_{1},\lambda_{2},\lambda_{3}) =0
\]
where $l$ is defined in equation (\ref{eqn:leqn}).

Before we find $\beta''$ we record a general formula which will be
 useful to us in many situations, if
\[ \zeta_{\lambda_{1},\lambda_{2}}=  E(x_{1} \lambda_{1} +x_{2} \lambda_{2}, x_{3} \lambda_{1} + x_{4} \lambda_{2}, x_{5} \lambda_{1} + x_{6} \lambda_{2})\]
for some coefficients $x_{i} \in \mathbb{C}$ then
\begin{equation}\label{eqn:helper}
\begin{split} & (\delta \zeta)_{\lambda_{1},\lambda_{2},\lambda_{3}}  \\ & = 2E(x_1 \lambda_{1}, (x_5 - x_6)\lambda_{2},x_4 \lambda_{3})
 + 2E(x_3 \lambda_{1}, (x_6 - x_5)\lambda_{2},x_2 \lambda_{3}) + 2E(x_1 \lambda_{1}, (x_4 - x_3)\lambda_{2},x_6 \lambda_{3}) \\
 & \quad + 2E(x_5 \lambda_{1}, (x_3 - x_4)\lambda_{2},x_2 \lambda_{3}) +2E(x_3 \lambda_{1}, (x_1 - x_2)\lambda_{2},x_6 \lambda_{3})
 + 2E(x_5 \lambda_{1}, (x_2 - x_1)\lambda_{2},x_4 \lambda_{3}).
\end{split}
\end{equation}

Using equation (\ref{eqn:helper}) it is easily seen that equation (\ref{kill_imag_part}) can be solved by
\begin{equation} \label{BetaDoublePrime}
\beta''_{\lambda_{1},\lambda_{2}} =
\frac{1}{16}\Big{(}E(i\lambda_{1},\lambda_{2},\lambda_{1})
-E(\lambda_{2},i(\lambda_{1}+\lambda_{2}),\lambda_{1})\Big{)}.
\end{equation}

Taking the real part of equation
(\ref{eqn:integral}) leaves us with
\begin{equation}\label{eqn:real_part}
(\delta \beta')_{\lambda_{1},\lambda_{2}, \lambda_{3}} +k(\lambda_{1},\lambda_{2}, \lambda_{3}) = 0
\end{equation}
or
\[
\beta'_{\lambda_{2},\lambda_{3}}
-\beta'_{\lambda_{1}+\lambda_{2},\lambda_{3}}
+\beta'_{\lambda_{1},\lambda_{2}+\lambda_{3}}
-\beta'_{\lambda_{1},\lambda_{2}}
+k(\lambda_{1},\lambda_{2}, \lambda_{3})
\in \mathbb{Z}.
\]
where $k$ is defined in equation (\ref{eqn:keqn}).
We can see that such a $\beta'$ always exists by the following argument.
  The short exact sequence
\[0 \to \mathbb{Z} \to \mathbb{R} \to U(1) \to 1
\]
leads to the following commutative diagram of short exact sequences where the vertical aarows are isomorphisms.
\[
\xymatrix{\ar @{} [dr] |{}
0 \ar[r] & H^{3}(\Lambda,\mathbb{Z}) \ar[d]^s \ar[r] &
  H^{3}(\Lambda,\mathbb{R}) \ar[d]^s \ar[r] & H^{3}(\Lambda,U(1)) \ar[d]^s \ar[r] & 1 \\
0 \ar[r] & \text{Alt}^{3}(\Lambda, \mathbb{Z}) \ar[r] &
\text{Alt}^{3}(\Lambda, \mathbb{R}) \ar[r] & \text{Alt}^{3}(\Lambda, U(1)) \ar[r] & 1}
\]
where the map $s$ on the groups involving $U(1)$ is induced from the first two vertical maps.

Observe that by linearity, $k \in Z^{3}(\Lambda, \mathbb{R})$.  Since
$k$ skew-symmetrizes to the integral element
\[s(k) = E \in \text{Alt}^{3}(\Lambda,\mathbb{Z}) \subset
  \text{Alt}^{3}(\Lambda,\mathbb{R}),
\] it corresponds to an integral cohomology class.  Hence $k$ is equivalent
 via an element of $B^{3}(\Lambda, \mathbb{R})$ in $Z^{3}(\Lambda, \mathbb{R})$
  to an element of $Z^{3}(\Lambda, \mathbb{Z})$.  Hence the image of $k$ in
   $Z^{3}(\Lambda, U(1))$ is trivializable by an element of $B^{3}(\Lambda, U(1))$.
     This precisely says that $\beta'$ exists.  For an
     explicit construction see the appendix in Section (\ref{appendix1}).  The role of $\beta'$ is similar to the role of the semi-character $\vartheta$ in Theorem \ref{lineAH}.  The solution for $\beta'$ corresponding to $E$ appears in equation (\ref{eqn:FinalBetaPrime}).
Finally, we can conclude that given an element $E \in A(\Lambda)$, the element $\Phi^{E}$ defined by
\begin{equation}\label{eqn:DefOfPhiE}
\Phi^{E}_{\lambda_{1},\lambda_{2}}(v)= \exp(\Theta^{E}_{\lambda_{1},\lambda_{2}}(v)) =\exp( H_{\lambda_{1},\lambda_{2}}(v) + \beta'_{\lambda_{1},\lambda_{2}}+i\beta''_{\lambda_{1},\lambda_{2}})
\end{equation}
of $Map(\Lambda \times \Lambda, \mathcal{O}^{\times}(V))$ lies in $Z^{2}(\Lambda, \mathcal{O}^{\times}(V))$ and the skew-symmetrization $s(\del \Theta^{E}) \in \text{Alt}^{3}(\Lambda, \mathbb{Z})$ of $ \del \Theta^{E} \in Z^{3}(\Lambda,\mathbb{Z})$ agrees with $E$.  In the following two definitions $E$ is an element of $A(\Lambda)$ which was defined in Definition \ref{Adefinitions} and $B$ is an element of $\wedge^{2}\overline{V}^{\vee}$.
\begin{dfn}\label{DefOfGerbeForE}
Let $\mathfrak{G}_{E}$ be the gerbe on $X$ corresponding via (\ref{eqn:ob2CyclesToGerbes}) to $\Phi^{E} \in Z^{2}(\Lambda, \mathcal{O}^{\times}(V))$ which was defined in equation (\ref{eqn:DefOfPhiE}).
\end{dfn}
\begin{dfn}
Let $\mathfrak{G}_{B}$ be the gerbe on $X$ corresponding via (\ref{eqn:ob2CyclesToGerbes}) to the constant cocycle $\Phi^{B} \in Z^{2}(\Lambda, \mathcal{O}^{\times}(V))$ given by
\[
\Phi^{B}_{\lambda_{1}, \lambda_{2}}(v) = \exp\Big{(}\frac{1}{2}B(\lambda_{1}, \lambda_{2})\Big{)}.
\]
\end{dfn}
We have a functor
\[F: [\wedge^{2}\overline{V}^{\vee} / \text{Alt}^{2}(\Lambda, \mathbb{Z}) ] \to \mathfrak{G}erbes(X)
\]
defined as follows, for every object $B \in \wedge^{2}\overline{V}^{\vee}$ we assign the gerbe $\mathfrak{G}_{B}$.  For any morphism $\mu \in \text{Alt}^{2}(\Lambda, \mathbb{Z})$ mapping $B_{1}$ to $B_{2}$, in other words $B_{2} - B_{1} = \mu^{H}$ we assign the isomorphism
\[F(\mu; B_{1}, B_{2}):\mathfrak{G}_{B_{1}} \to \mathfrak{G}_{B_{2}}
\]
determined in the sense of (\ref{eqn:mor2CyclesToGerbes}) by the cocycle in $C^{1}(\Lambda, \mathcal{O}^{\times}(V))$ given by the same formula as the cocycle in equation (\ref{eqn:IsomInGerbes}).  The functor $F$ here is considered as a functor between $2-$categories where the $2-$morphisms on the left hand side are only the identities.  Similarly,

\begin{thm}
There is a functor (of 2-categories)
\[\bold{ah} : [\wedge^{2}\overline{V}^{\vee} / \textup{Alt}^{2}(\Lambda, \mathbb{Z}) ] \times A(\Lambda) \to \mathfrak{G}erbes(X).
\]
Such that for any object $(B, E)$, the gerbe $\bold{ah}(B, E)$ has topological class corresponding to $E$ and the cohomology class $[\bold{ah}(B, 0)]$ is the image of $B \in H^{2}(X, \mathcal{O})$ inside $H^{2}(X, \mathcal{O}^{\times})$.  The induced map
\[(\wedge^{2}\overline{V}^{\vee} / \textup{Alt}^{2}(\Lambda, \mathbb{Z})^{H})  \times A(\Lambda) \to H^{2}(X, \mathcal{O}^{\times})\] is an isomorphism of groups where the group structure on the left hand side is induced from the groups $\wedge^{2}\overline{V}^{\vee}$ and $\text{Alt}^{3}(\Lambda, \mathbb{Z})$ and on the right hand side the group structure comes from the monoidal structure on gerbes.
\end{thm}
{\bf Proof.}
On the level of objects, the map is
\[(B,E) \mapsto F(B) \otimes \mathfrak{G}_{E}.
\]
On the level of morphisms, we have
\[\text{Hom}((B_{1},E),(B_{2},E)) \to \text{Hom}(F(B_{1}) \otimes \mathfrak{G}_{E},F(B_{2}) \otimes \mathfrak{G}_{E})
\]
given by
\[\mu \mapsto F(\mu; B_{1}, B_{2})\otimes \text{id}.
\]
The required properties have all been proven above.

\ \hfill $\Box$

\noindent
As groups we have
 \[
 \Big{(}\wedge^{2}\overline{V}^{\vee} / \textup{Alt}^{2}(\Lambda, \mathbb{Z})^{H}\Big{)} = \Big{(}\text{Alt}^{2}(\Lambda, \mathbb{R})/(\text{Alt}^{2}(\Lambda, \mathbb{R})^{(1,1)}+\text{Alt}^{2}(\Lambda,\mathbb{Z})) \Big{)}.\]
\begin{rmk}\label{rmk:important}
Notice that for every gerbe we have constructed a cocycle representative based on the choice of $B$ and $E$, but there is no unique choice because of the kernel of the Hodge projection acting on $B$.  We name these cocycle representatives $\Phi^{(B,E)} \in Z^{2}(\Lambda, \mathcal{O}^{\times}(V))$.
\end{rmk}
\begin{dfn}\label{dfn:Reps}
The cocycle $\Phi^{(B,E)} = \Phi^{B}\Phi^{E}$ is defined by
\[\Phi^{(B,E)}_{\lambda_{1}, \lambda_{2}}(v)=\exp( \frac{1}{2}B(\lambda_{1},\lambda_{2})+ H_{\lambda_{1},\lambda_{2}}(v) + \beta'_{\lambda_{1},\lambda_{2}}+i\beta''_{\lambda_{1},\lambda_{2}})
\]
where $H$ is defined in (\ref{eqn:Heqn}), $\beta'$ is defined in (\ref{eqn:FinalBetaPrime}), and $\beta''$ is defined in (\ref{BetaDoublePrime}).
\end{dfn}
\section{A Groupoid in Sheaves of Sets}\label{groupoid}
The cocycle we have constructed $\Phi=\Phi^{(B,E)}: \Lambda \times \Lambda \to \mathcal{O}^{\times}(V)$ is normalized in the sense that $\Phi_{\lambda_{1},0} = \Phi_{0,\lambda_{2}} = 1$ for all
$\lambda_{1}, \lambda_{2} \in \Lambda$.  It is easy to see that any such cocycle defines a groupoid with structure maps

\begin{equation}\label{GroupoidForGerbe}
\xymatrix{**[l] V \times \Lambda \times \mathbb{C}^{\times}  \ar@/_2pc/[r]_-{s} \ar@/^2pc/[r]^-{t} & **[r] \ar@/^0pc/[l]_-{e}  V}
\end{equation}

\[\xymatrix{{(V\times \Lambda \times \mathbb{C}^{\times})}_{t} \times_{V} {{{}_{s}}(V\times \Lambda \times \mathbb{C}^{\times})} \ar@/^0pc/[r]^-{m} & V}
\]
and
\[\xymatrix{(V \times \Lambda \times \mathbb{C}^{\times}) \ar@/^0pc/[r]^-{\iota} & (V \times \Lambda \times \mathbb{C}^{\times})  }
\]
Here $s$ is given by
\[s(v,\lambda,z)=v
\]
$t$ is given by
\[t(v,\lambda,z)=v+\lambda
\]
$e$ is given by
\[e(v) = (v,0,1)
\]
the inverse map is
\[\iota_{\Phi}(v,\lambda,z) =
(v+\lambda,-\lambda,z^{-1}\Phi_{\lambda,-\lambda}^{-1}(v))
\]
and the multiplication map is given by
\[m_{\Phi}\Big{(}(v_{1},\lambda_{1},z_{1}),(v_{2} = v_{1}+ \lambda_{1},\lambda_{2},z_{2})\Big{)} = (v_{1},\lambda_{1}+\lambda_{2}, z_{1}z_{2} \Phi_{\lambda_{1}, \lambda_{2}}(v_{1})).
\]

Suppose we have a groupoid with objects $G_{0}$ and morphisms $G_{1}$ in the category of complex analytic spaces and maps from the source and target to a complex analytic space $Y$ such that the following diagram commutes.  
\[\xymatrix{&**[l] G_{1}  \ar[ddr] \ar@/_2pc/[rr]_-{s} \ar@/^2pc/[rr]^-{t} & & **[r] \ar@/^0pc/[ll]_-{e} \ar[ddl] G_{0} \\ \\
& & Y}
\]
Now the sheaves of sections of the maps to $Y$ form a groupoid in the category of sheaves of sets on $Y$.  The source, target, multiplication, and identity maps are all the obvious induced ones on sheaves of sections.  Let $X = V/\Lambda$ be a complex torus and $\mathcal{V}^{p}_{X}$ denote the sheaf on $X$ given by sections of the map $p:V \to X$.  Let $\Lambda_{X}$ be the constant sheaf of groups on $X$ with fiber $\Lambda$.  Starting with our groupoid in equation (\ref{GroupoidForGerbe}) we get the following groupoid in the category of sheaves of sets on $X$
\begin{equation}\label{eqn:weird_thing}
\xymatrix{**[l] \mathcal{V}^{p}_{X} \times \Lambda_{X} \times \mathcal{O}_{X}^{\times}  \ar@/_2pc/[r]_-{s} \ar@/^2pc/[r]^-{t} &  **[r] \ar@/^0pc/[l]_-{e}  \mathcal{V}^{p}_{X}}.
\end{equation}

Recall that the action of a groupoid with objects $G_{0}$ and morphisms $G_{1}$ on a set $T$ is a map $a_{0}:T \to G_{0}$ and a map $a_{1}:{G_{1}}_{s} \times_{G_{0}} {}_{a_{0}}T \to T$ compatible with the composition map $m:G_{1 s} \times_{G_{0}} {}_{t}G_{1} \to G_{1}$.  It is simply transitive whenever given $t_{1}, t_{2} \in T$ there is a unique $\gamma \in G_{1}$ such that $s(\gamma) = a_{0}(t_{1})$ and $a_{1}(\gamma, t_{1}) = t_{2}$.  Given a groupoid in sheaves of sets on $X$, a torsor for this groupoid is a sheaf of sets on $X$ along with an action of the groupoid in sheaves which is locally simply transitive.  Given a groupoid in the category of sheaves of sets over $X$, the category of torsors on any open set forms a groupoid.

Denote the resulting stack on $X$ of torsors for (\ref{eqn:weird_thing}) by $\mathfrak{G}(\Phi)$.  We claim that $\mathfrak{G}({\Phi})$ is a gerbe and represents the class of $\Phi$ in
cohomology $H^{2}(\Lambda,\mathcal{O}(V)^{\times}) \cong H^{2}(X,\mathcal{O}^{\times})$.  In order to see this we can write down an explicit isomorphism from the gerbe $\mathfrak{G}_{\Phi}$ assoiciated (\ref{eqn:2CyclesToGerbes}) to $\Phi$ in the second appendix, to the torsors $\mathfrak{G}({\Phi})$ for (\ref{eqn:weird_thing}).  Given an open set $U \subset X=V/\Lambda$, $\mathfrak{G}_{\Phi}$ associates the groupoid whose objects are pairs of an $\mathcal{O}^{\times}$-torsor $\Xi$ on the principal $\Lambda$ bundle $p^{-1}(U)\to U$ and and isomorphisms
\[G_{\lambda}:   \lambda_{*} \Xi \to \Xi
\]
satisfying an obvious compatibility with $\Phi$.  We can actually see that $\Xi \to U$ is a torsor over the groupoid in sheaves (\ref{eqn:weird_thing}) restricted to $U$.  There is an obvious map $a_{0}$ which takes sections of $\Xi$ over $U$ to sections of $p^{-1}(U)$ over $U$.  Now if $q$ is a local section of $\Xi$ over $U$ covering a local section
$\varsigma$ of $p^{-1}(U)$ over $U$, in other words $a_{0}(q) = \varsigma$, then the triple $(\varsigma, \lambda, f)$ of sections over $U$ acts on $q$ by taking it to
\[a_{1}((\varsigma, \lambda, f);q)=f G_{\lambda}(\lambda_{*}  q).\]  In order to check that this is indeed an action, observe that for any section $q$ of $\Xi \to X$ we have
\begin{equation}
\begin{split}
& a_{1}\Big{(}(\varsigma,\lambda_{1},f_{1});a_{1}(\varsigma,\lambda_{2},f_{2};q)\Big{)} \\
&= a_{1}\Big{(}(\varsigma,\lambda_{1},f_{1}) ; f_{2}G_{\lambda_{2}} (\lambda_{2*}q)\Big{)} \\
& = f_{1} f_{2} G_{\lambda_{1}} (\lambda_{1*}G_{\lambda_{2}})(\lambda_{1}+\lambda_{2})_{*} q \\
& =  f_{1} f_{2} G_{\lambda_{1}+\lambda_{2}} (\Phi_{\lambda_{1},\lambda_{2}}(\varsigma))(\lambda_{1}+\lambda_{2})_{*} q \\
& = a_{1}\Big{(}\varsigma,\lambda_{1}+\lambda_{2},
f_{1}f_{2}(\Phi_{\lambda_{1},\lambda_{2}}(\varsigma));q\Big{)} \\
& = a_{1}\Big{(}m_{\Phi}((\varsigma,\lambda_{1},f_{1}),(\varsigma+\lambda_{1},\lambda_{2},f_{2}));q\Big{)}
\end{split}
\end{equation}
as needed.  We conclude that $p_{*}\Xi$ is a torsor for (\ref{eqn:weird_thing}) for every $\Xi$ in $\mathfrak{G}_{\Phi}(U)$.  This sets up an isomorphism $\mathfrak{G}_{\Phi} \to \mathfrak{G}(\Phi)$.
\section{The Interaction With the Group Structure}\label{mult}
Let $S$ be a complex analytic space.  Let $F$ be a contravariant
functor from the category of complex analytic spaces $X$ over $S$
with section $x$ to abelian groups.  Consider any collection of
complex analytic spaces $X_{i}$ over $S$ with sections $x_{i}$.  Any
fiber product of these over $S$ will be considered a complex analytic
space over $S$ with the obvious section.  We have the
inclusion morphisms
\[
\sigma_{i}^{n} : X_{0} \times_{S} \cdots \times_{S} \widehat{X_{i}} \times_{S}
\cdots \times_{S} X_{n} \to  X_{0} \times_{S} \cdots \times_{S} X_{n}
\]
defined by using the section $x_{i}:S \to X_{i}$ in the $i-$th location.  Applying
the contravariant functor $F$, we get group homomorphisms
 \[F(\sigma,n)=\prod_{i=0}^{n} F(\sigma_{i}^{n}):F( X_{0} \times_{S} \cdots \times_{S} X_{n}) \to
 F(X_{0} \times_{S} \cdots \times_{S} \widehat{X_{i}} \times_{S} \cdots \times_{S} X_{n} ).
\]
Following \cite{Mu1970}, $F$ is said to be {\it of order $n$} if the group
homomorphisms $F(\sigma, n)$ is injective.   For
$n=0$ such a $F$ will be called constant, for $n=1$ such a $F$ will be
called linear, for $n=2$ quadratic, and for $n=3$ cubic.  It is easy to see that if $F$ is of order $n$ then it is also of order $m$ for all $m>n$.

Let us say that a map
$f:X \to S$ {\it satisfies condition $C_{p}$} if for every point $s \in S$ there is an arbitrarily small contractible
neighborhood $U$
containing $s$ such that (i) the inclusion map of the fiber $f^{-1}(s) \hookrightarrow f^{-1}(U)$
induces an isomorphism
\begin{equation}\label{eqn:TopConst1}
H^{p}(f^{-1}(U),\mathbb{Z}) \cong H^{p}(f^{-1}(s),\mathbb{Z})
\end{equation}
and (ii) the natural map of restricting to the fibers
\begin{equation}\label{eqn:TopConst2}
H^{p}(f^{-1}(U),\mathbb{Z}) \to H^{0}(U, R^{p}f_{*} \mathbb{Z})
\end{equation}
is an isomorphism.

\begin{lem}\label{lem:TopConstraints}
If $f:X \to S$ satisfies conditions $C_{p}$ and $C_{p+1}$ then an element of
$H^{0}(S, R^{p}f_{*} \mathcal{O}^{\times})$ is trivial if and only if it is trivial when restricted to each fiber.
\end{lem}
{\bf Proof.}
We will use the exact sequence
\[R^{p}f_{*} \mathbb{Z} \to R^{p}f_{*} \mathcal{O} \to R^{p}f_{*} \mathcal{O}^{\times} \to R^{p+1}f_{*} \mathbb{Z}
\]
Choose any $\sigma \in H^{0}(S, R^{p}f_{*} \mathcal{O}^{\times})$ and assume that it is trivial
when restricted to each fiber.  That is to say that its image $\sigma_{s} \in
H^{p}(f^{-1}(s), \mathcal{O}^{\times})$ is trivial for each $s \in S$.  If we fix $s \in S$ arbitrary, it suffices  to produce a neighborhood $U$ of $s$ on which $\sigma$ trivializes.  Chose $U'$ contractible containing $s$ such that both the vertical and the horizontal maps on the right hand side of the commutative diagram
 \[
\xymatrix{\ar @{} [dr] |{}
  &
  H^{0}(U',R^{p}f_{*}\mathcal{O}^{\times} )  \ar[d] \ar[r] & H^{0}(U',R^{p+1}f_{*}\mathbb{Z}) \ar[d] &  H^{p+1}(f^{-1}(U'),\mathbb{Z}) \ar[l]&\\
  &
 H^{p}(f^{-1}(s), \mathcal{O}^{\times})  \ar[r] & H^{p+1}(f^{-1}(s),\mathbb{Z})  & &}
\]
 are isomorphisms.  This diagram then shows that we can pick a pre-image $\tilde{\sigma} \in H^{0}(U',R^{p}f_{*}\mathcal{O} )$ to $\sigma|_{U'}$.  The image $\tilde{\sigma}_{s}$ of  $\tilde{\sigma}$ in each fiber $H^{p}(f^{-1}(s), \mathcal{O})$ actually comes from an element of $H^{p}(f^{-1}(s), \mathbb{Z})$.  Chose a contractible neighborhood $U$ of $s$ inside $U'$ such that
both the vertical and the horizontal maps on the left hand side of the diagram
\[
\xymatrix{\ar @{} [dr] |{}
  & H^{p}(f^{-1}(U),\mathbb{Z}) \ar[r] &
  H^{0}(U,R^{p}f_{*}\mathbb{Z} )  \ar[d] \ar[r] & H^{0}(U,R^{p+1}f_{*}\mathcal{O}) \ar[d] &  \\
  & &
 H^{p}(f^{-1}(s), \mathbb{Z})  \ar[r] & H^{p+1}(f^{-1}(s),\mathcal{O})  & }
\]
are isomorphisms.  This shows finally that $\tilde{\sigma}|_{U}$ comes from $H^{0}(U,R^{p}f_{*}\mathbb{Z} )$ and hence that $\sigma$ is trivial when restricted to $U$.  Hence $\sigma$ is trivial.

\ \hfill $\Box$

\begin{lem}\label{lem:ProperTop}Let $f:X \to S$ be any proper analytic map of complex analytic spaces.  Then an element of
$H^{0}(S, R^{p}f_{*} \mathcal{O}^{\times})$ is trivial if and only if it is trivial when restricted to each fiber.
\end{lem}

{\bf Proof.}
This follows immediately from two theorems on {\it constructible sheaves} of abelian groups with respect to an analytic Whitney stratification of a complex analytic space.  The first (Theorem 4.1.5 (i) (b)
\cite{Di2004}) says that the sheaves $R^{p}f_{*} \mathbb{Z}$ (the push-forwards of the constructible sheaf $\mathbb{Z}$ under a proper analytic map) are all constructible.  The second (Theorem 4.1.9) \cite{Di2004}) says that for any constructible sheaf $\mathcal{S} \to M$ on a complex analytic space $M$, and any point $m \in M$, there exists arbitrarily small contractible neighborhoods $U$ with $m \in U \subset M$ and such that (i) $H^{p}(U,\mathcal{S})=0$ for all $p>0$ and (ii) the inclusion map of $m$ into $U$ induces an isomorphism $H^{0}(U,\mathcal{S}) \cong \mathcal{S}_{m}$.  Now by (i) the Leray spectral sequence for $H^{p}(f^{-1}(U), \mathbb{Z})$ collapses and we get the isomorphisms in (\ref{eqn:TopConst2}).  Finally, we use the fact that $(R^{p}f_{*} \mathbb{Z})_{s} \cong H^{p}(f^{-1}(s), \mathbb{Z})$ for any $s \in S$ and therefore (ii) imply the isomorphism in (\ref{eqn:TopConst1}).  Therefore the topological constraints $C_{p}$ and $C_{p+1}$ are satisfied by $f$ and we are done by Lemma \ref{lem:TopConstraints}.

\ \hfill $\Box$

\begin{thm}\label{thm:CubicFunctor}The functor $F(X) = H^{2}(X,\mathcal{O}^{\times})$ is cubic when restricted to the full subcategory of complex analytic spaces $X$, with a proper map $f: X \to S$ with section satisyfing $f_{*} \mathcal{O} = \mathcal{O}$.  On this same category, $H^{1}(X,\mathcal{O}^\times)$ is quadratic.
\end{thm}
Notice that condition $f_{*} \mathcal{O} = \mathcal{O}$ follows if all fibers are reducible.    We remark that this theorem was inspired by a similar result in a paper \cite{Ho1972} of R. Hoobler.  There he proves in an algebraic setting (where gerbes are topologically trivial) that the analogous functor is quadratic.  When our gerbes are topologically trivial along the fibers, our results agree with his.  Its possible that our results could be use to extend the validity of his results to other cases in the algebraic setting over $\mathbb{C}$.  Before proving the theorem, we will need a lemma proving it in the special case that $S$ is a point.

\begin{lem}\label{lem:NotMe}
It follows from the K\"unneth decompositions for $\mathbb{Z}$ (see \cite{Di2004}) and $\mathcal{O}$ (see the book of Demailly \cite{De} in the section on Grauert's direct image theorem) that for compact complex analytic spaces in the absolute setting ($S$ is a point) the functor \[X \mapsto H^{p}(X,\mathbb{Z})/(tors(H^{p}(X,\mathbb{Z}))\] is cubic for $p=2,3$ and quadratic for $p=2$.  Similarly, the functor \[X \mapsto H^{q}(X, \mathcal{O})\] is both quadratic and cubic for $q=1,2$.  Therefore the functor $X \mapsto H^{1}(X,\mathcal{O}^{\times})$ is quadratic and the functor $X \mapsto H^{2}(X,\mathcal{O}^{\times})$ is cubic.
\end{lem}

\ \hfill $\Box$

{\bf Proof of Theorem \ref{thm:CubicFunctor}.}
Let  $f: Y \to S$ be a morphism of complex analytic spaces satisfying 
$f_{*}\mathcal{O} = \mathcal{O}$ such that $f$ has a section.    For each $i$ the pullback maps 
\[f^{*}: H^{i}(S,\mathcal{O}^{\times}) \to H^{i}(Y,\mathcal{O}^{\times})
\]
have a left inverse provided by the pullback with respect to the section.  In particular they are injective maps of abelian groups.  A sector of the $E_{2}$ term of the spectral sequence looks as follows where we have identified $f_{*}\mathcal{O}_{X}^{\times}$ with $\mathcal{O}_{S}^{\times}$.
\def\csg#1{\save[].[ddrr]!C*+<2pc,1pc>[F-,]\frm{}\restore}
\[
\xymatrix@=1.5pc{
q=2 & \csg1  H^{0}(S,R^{2}f_{*}\mathcal{O}^{\times})
&  H^{1}(S,R^{2}f_{*}\mathcal{O}^{\times}) &
H^{2}(S,R^{2}f_{*}\mathcal{O}^{\times}) 
\\
q= 1 & H^{0}(S,R^{1}f_{*}\mathcal{O}^{\times})
&  H^{1}(S,R^{1}f_{*}\mathcal{O}^{\times}) &
H^{2}(S,R^{1}f_{*}\mathcal{O}^{\times})  \\
q=0 & H^{0}(S,\mathcal{O}^{\times}) & H^{1}(S,\mathcal{O}^{\times}) &H^{2}(S,\mathcal{O}^{\times}) \\
 & p= 0 & p= 1 & p= 2 
}
\]
The injectivity of the pullback maps imply that all differentials landing in the bottom row must be zero and hence the bottom row does not change from page to page.  The same sector on the $E_{\infty}$ page looks like 
\def\csg#1{\save[].[ddrr]!C*+<4.5pc,1pc>[F-,]\frm{}\restore}
\[
\xymatrix@=1.5pc{
q=2 & \csg1  \text{ker}(d_2)
& \star &
\star
\\
q= 1 & H^{0}(S,R^{1}f_{*}\mathcal{O}^{\times})
&  H^{1}(S,R^{1}f_{*}\mathcal{O}^{\times}) &
H^{2}(S,R^{1}f_{*}\mathcal{O}^{\times})/\text{im}(d_2)   \\
q=0 & H^{0}(S,\mathcal{O}^{\times}) & H^{1}(S,\mathcal{O}^{\times}) & H^{2}(S,\mathcal{O}^{\times})  \\
 & p= 0 & p= 1 & p= 2 
}
\]
We know that there is a filtration $F^{i}H^{2}(Y,\mathcal{O}^{\times})$ such that the $p+q=2$ terms are identified with the associated graded groups.  Therefore  $F^{0}H^{2}(Y,\mathcal{O}^{\times}) =H^{2}(S,\mathcal{O}^{\times})$,  $F^{2}H^{2}(Y,\mathcal{O}^{\times}) =H^{2}(Y,\mathcal{O}^{\times})$, 
\[F^{2}H^{2}(Y,\mathcal{O}^{\times}) /F^{1} H^{2}(Y,\mathcal{O}^{\times}) = \text{ker} [d_{2}:H^{0}(S,R^{2}f_{*}\mathcal{O}^{\times}) \to H^{2}(S, R^{1}f_{*}\mathcal{O}^{\times})]
\]
and there is an exact sequence 
\[1 \to H^{2}(S,\mathcal{O}^{\times})  \to F^{1}H^{2}(Y,\mathcal{O}^{\times}) \to H^{1}(S,R^{1}f_{*}\mathcal{O}^{\times}) \to 1.
\]
The section gives a map $H^{2}(Y,\mathcal{O}^{\times}) \to H^{2}(S,\mathcal{O}^{\times})$ and by restriction to the middle term, can be used to split this sequence.  Thus we have
\[ H^{2}(S,\mathcal{O}^{\times}) \oplus H^{1}(S, R^{1}f_{*} \mathcal{O}^{\times}) \subset H^{2}(Y,\mathcal{O}^{\times})
\]
where the quotient
\begin{equation} \label{eqn:DefOfQ}
Q(Y,f) = H^{2}(Y,\mathcal{O}^{\times})/(H^{2}(S,\mathcal{O}^{\times}) \oplus H^{1}(S, R^{1}f_{*} \mathcal{O}^{\times}))
\end{equation}
lies in $H^{0}(S,R^{2}f_{*} \mathcal{O}^{\times})$.
Let us start with a class in
\[H^{2}(X_{1} \times_{S} X_{2} \times_{S} X_{3} \times_{S} X_{4},\mathcal{O}^{\times})\]
which is trivial when restricted to each of the 'slices' $X_{i} \times_{S} X_{j} \times_{S} X_{k}$ for $1 \leq i < j < k \leq 4$.
We need to show that this class is trivial.  Let
\[G \in H^{0}(S,R^{2}p^{1,2,3,4}_{*} \mathcal{O}^{\times})
\]
be the restriction of the class to the fibers.  Clearly, the restriction of $G$ to each slice of each fiber  of
\[p^{1,2,3,4}: X_{1} \times_{S} X_{2} \times_{S} X_{3} \times_{S} X_{4} \to S
\] is trivial.   Therefore by lemma \ref{lem:ProperTop} and the fact that our theorem already holds for the case that $S$ is a point (see lemma \ref{lem:NotMe}), we see that $G$ is trivial.  So we have concluded so far that the functor $X \mapsto Q(X,f)$ is cubic.  It therefore will suffice to prove that the functor $X \mapsto H^{2}(S,\mathcal{O}^{\times}) \oplus H^{1}(S, R^{1}f_{*} \mathcal{O}^{\times})$ is cubic.  The factor $X \mapsto H^{2}(S,\mathcal{O}^{\times})$ is constant and hence cubic.  It will therefore suffice to show that the functor \[X \mapsto H^{1}(S, R^{1}f_{*} \mathcal{O}^{\times})\] is quadratic, since a quadratic functor is cubic.  Therefore we consider now the projection
\[p^{1,2,3}: X_{1} \times_{S} X_{2} \times_{S} X_{3} \to S.\] Suppose
therefore that we have an element \[G \in H^{1}(S,R^{1}p^{1,2,3}_{*} \mathcal{O}^{\times})\] which
is trivial when restricted to each of the spaces $X_{i} \times_{S} X_{j}$ for $1 \leq i < j \leq 3$.  We need to show that $G$ is trivial.

We will use the notation \[p^{i}:X_{i} \to S,\]
\[p^{i,j}:X_{i} \times_{S} X_{j} \to S,\]
\[p^{i,j}_{k}:X_{i} \times_{S} X_{j} \to X_{k},\]
\[p^{i,j,k}_{n}:X_{i} \times_{S} X_{j} \times_{S} X_{k} \to X_{n}\] and
\[p^{i,j,k}_{n,m}:X_{i} \times_{S} X_{j} \times_{S} X_{k} \to X_{n} \times_{S} X_{m}\] to denote the projections.
We can pick and open cover $\{ U_{\alpha} \}$ of $S$ such that we can represent $G$ by a \v{C}ech cocycle, in other words by a collection of relative line bundles on $({{p^{1,2,3}}})^{-1}(U_{\alpha})$.  In other words, we pick \[L_{\alpha,\beta} \in H^{0}(U_{\alpha,\beta}, R^{1}p^{1,2,3}_{*} \mathcal{O}^{\times})\] satisfying
\[L_{\alpha,\beta} \otimes L_{\beta,\gamma} = L_{\alpha,\gamma}
\] which represents $G$.  For each $i,j$ such that $1 \leq i < j  \leq 3$ we know that there exists \[M_{\alpha}^{i,j} \in H^{0}(U_{\alpha}, R^{1}p^{i,j}_{*} \mathcal{O}^{\times})\] such that the restrictions \[L^{i,j}_{\alpha,\beta} = L_{\alpha,\beta}|_{X_{i} \times_{S} X_{j}}\] of $L_{\alpha,\beta}$ to $X_{i} \times_{S} X_{j}$ satisfy
\begin{equation}\label{eqn:BoundOnSlices}L^{i,j}_{\alpha,\beta} = ((M_{\alpha}^{i,j})\otimes {(M_{\beta}^{i,j})}^{-1})|_{U_{\alpha,\beta}}.
\end{equation}
 Define
 \[N^{1}_{\alpha} = M^{1,2}_{\alpha} |_{X_{1}} \in
 H^{0}(U_{\alpha}, R^{1}p^{1}_{*} \mathcal{O}^{\times}),
 \]
 \[P^{1}_{\alpha} = M^{1,3}_{\alpha} |_{X_{1}} \in
 H^{0}(U_{\alpha}, R^{1}p^{1}_{*} \mathcal{O}^{\times}),
 \]
 \[N^{2}_{\alpha} = M^{1,2}_{\alpha} |_{X_{2}} \in
 H^{0}(U_{\alpha}, R^{1}p^{2}_{*} \mathcal{O}^{\times}),
 \]
 \[Q^{2}_{\alpha} = M^{2,3}_{\alpha} |_{X_{2}} \in
 H^{0}(U_{\alpha}, R^{1}p^{2}_{*} \mathcal{O}^{\times}),
 \]
  \[P^{3}_{\alpha} = M^{1,3}_{\alpha} |_{X_{3}} \in
  H^{0}(U_{\alpha}, R^{1}p^{3}_{*} \mathcal{O}^{\times})
 \]
 and
 \[Q^{3}_{\alpha} = M^{2,3}_{\alpha} |_{X_{3}} \in
 H^{0}(U_{\alpha}, R^{1}p^{3}_{*} \mathcal{O}^{\times}).
 \]
Then we compute
\[((N^{1}_{\alpha})\otimes (P^{1}_{\alpha})^{-1})|_{U_{\alpha,\beta}} \otimes ((N^{1}_{\beta})^{-1} \otimes (P^{1}_{\beta}))_{U_{\alpha,\beta}} = (L^{1,2}_{\alpha,\beta})|_{X_{1}} \otimes (L^{1,3}_{\alpha,\beta})^{-1}|_{X_{1}} = \mathcal{O}_{({p^{1})}^{-1}(U_{\alpha,\beta})},
\]
\[((N^{2}_{\alpha})\otimes (Q^{2}_{\alpha})^{-1})|_{U_{\alpha,\beta}} \otimes ((N^{2}_{\beta})^{-1} \otimes (Q^{2}_{\beta}))|_{U_{\alpha,\beta}} = (L^{1,2}_{\alpha,\beta})|_{X_{2}} \otimes (L^{2,3}_{\alpha,\beta})^{-1}|_{X_{2}} = \mathcal{O}_{({p^{2})}^{-1}(U_{\alpha,\beta})},
\]
and
\[((P^{3}_{\alpha})\otimes (Q^{3}_{\alpha})^{-1})|_{U_{\alpha,\beta}} \otimes ((P^{3}_{\beta})^{-1} \otimes (Q^{3}_{\beta}))|_{U_{\alpha,\beta}} = (L^{1,3}_{\alpha,\beta})|_{X_{3}}\otimes (L^{2,3}_{\alpha,\beta})^{-1}|_{X_{3}} = \mathcal{O}_{({p^{3})}^{-1}(U_{\alpha,\beta})}.
\]
Define
\[T^{1} \in H^{0}(S, R^{1}p^{1}_{*}\mathcal{O}^{\times})
\]
to be the element that restricts on every $U_{\alpha}$ to $(N_{\alpha}^{1}) \otimes (P_{\alpha}^{1})^{-1}$,
define
\[T^{2} \in H^{0}(S, R^{1}p^{2}_{*}\mathcal{O}^{\times})
\]
to be the element that restricts on every $U_{\alpha}$ to $(Q_{\alpha}^{2}) \otimes (N_{\alpha}^{2})^{-1}$ and define
\[T^{3} \in H^{0}(S, R^{1}p^{3}_{*}\mathcal{O}^{\times})
\]
to be the element that restricts on every $U_{\alpha}$ to $(P_{\alpha}^{3}) \otimes (Q_{\alpha}^{3})^{-1}$.

Then we have
\[M^{1,2}_{\alpha}|_{X_{1}} = M^{1,3}_{\alpha}|_{X_{1}} \otimes T^{1}|_{U_{\alpha}},
\]
\[M^{2,3}_{\alpha}|_{X_{2}} = M^{1,2}_{\alpha}|_{X_{2}} \otimes T^{2}|_{U_{\alpha}}
\]
and
\[M^{1,3}_{\alpha}|_{X_{3}} = M^{2,3}_{\alpha}|_{X_{3}} \otimes T^{3}|_{U_{\alpha}}.
\]
Define
\[W^{1,3}_{\alpha} = M^{1,3}_{\alpha} \otimes {(p^{1,3}_{1})}^{*} T^{1}|_{U_{\alpha}},
\]
\[W^{1,2}_{\alpha} = M^{1,2}_{\alpha} \otimes {(p^{1,2}_{2})}^{*} T^{2}|_{U_{\alpha}}
\]
and
\[W^{2,3}_{\alpha} = M^{2,3}_{\alpha} \otimes {(p^{2,3}_{3})}^{*} T^{3}|_{U_{\alpha}}.
\]
Then we have (just as in equation (\ref{eqn:BoundOnSlices}))
\begin{equation} \label{eqn:JustLike}
L^{i,j}_{\alpha,\beta} = ((W_{\alpha}^{i,j})\otimes {(W_{\beta}^{i,j})}^{-1})|_{U_{\alpha,\beta}}.
\end{equation}
Also the $W$ terms agree on their common factors and so we can define
\[W^{1}_{\alpha}:= W^{1,2}_{\alpha}|_{X_{1}} = W^{1,3}_{\alpha}|_{X_{1}},
\]
\[W^{2}_{\alpha}:= W^{1,2}_{\alpha}|_{X_{2}} = W^{2,3}_{\alpha}|_{X_{2}},
\]
and
\[W^{3}_{\alpha}:= W^{2,3}_{\alpha}|_{X_{3}} = W^{1,3}_{\alpha}|_{X_{3}}.
\]
Finally, define $M_{\alpha} \in H^{0}(U_{\alpha}, R^{1}p^{1,2,3}_{*} \mathcal{O}^{\times})$ by
\begin{equation}
\begin{split}
M_{\alpha}  = &({(p^{1,2,3}_{1,2})}^{*}W^{1,2}_{\alpha}) \otimes ({(p^{1,2,3}_{1,3})}^{*}W^{1,3}_{\alpha}) \otimes ({(p^{1,2,3}_{2,3})}^{*}W^{2,3}_{\alpha})\otimes \\
&
({(p^{1,2,3}_{1})}^{*}W^{1}_{\alpha})^{-1} \otimes ({(p^{1,2,3}_{2})}^{*}W^{2}_{\alpha})^{-1} \otimes ({(p^{1,2,3}_{3})}^{*}W^{3}_{\alpha})^{-1}.
\end{split}
\end{equation}
It is easily checked that
\[M_{\alpha}|_{X_{i} \times_{S} X_{j}} = W^{i,j,}_{\alpha}
.\]
We want to prove that
\begin{equation}\label{eqn:WhatWeWant}
M_{\alpha} \otimes M_{\beta}^{-1} = L_{\alpha, \beta}.
\end{equation}
First of all using Lemma \ref{lem:TopConstraints} it suffices to prove equation (\ref{eqn:WhatWeWant}) after restricting to each fiber
\[(p^{1,2,3})^{-1}(s) = (p^{1})^{-1}(s) \times (p^{2})^{-1}(s) \times (p^{3})^{-1}(s).
\]
In fact, using lemma \ref{lem:NotMe} it suffices to prove equation (\ref{eqn:WhatWeWant}) after restricting to each slice of each fiber, the slices being $(p^{i})^{-1}(s) \times (p^{j})^{-1}(s)$ for $1 \leq i <j \leq 3$.  Since the restriction of $M_{\alpha}$ to $(p^{i})^{-1}(s) \times (p^{j})^{-1}(s)$ is nothing but $W^{i,j}_{\alpha}|_{s}$, we are done in light of equation (\ref{eqn:JustLike}) which implies that
\[(W^{i,j}_{\alpha}|_{s}) \otimes {(W^{i,j}_{\beta}|_{s})}^{-1}  = L^{i,j}_{\alpha, \beta}|_{s}.
\]
Therefore $G$ is trivial and so we are done.  The proof that $X \mapsto H^{1}(X,\mathcal{O}^{\times})$ is quadratic uses precisely the same method as the proof that $X \mapsto Q(X,f)$ is cubic.

\ \hfill $\Box$

From this proof we get the following corollary.

\begin{cor}\label{cor:filtration}
Let $F$ be the contravariant functor from complex analytic spaces $f:X \to S$ with section to
abelian groups given by $F(X) = H^{2}(X,\mathcal{O}^{\times})$.   Then there is
short exact sequence of functors
\[1 \to G \to F \to G_{3} \to 1
\]
and a decomposition
\[G= G_{0} \oplus G_{2}
\]
such that
$G_{0}$ is constant, $G_{2}$ is quadratic and $G_{3}$ is cubic.  Furthermore, $G_{3}$ fits into an exact sequence
\[1 \to K_{2} \to G_{3} \to K_{3} \to 1
\]
where $K_{2}$ is quadratic and $K_{3}$ is cubic.
We also observe
that $G_{2}$ admits a filtration
\[1 \to H_{1} \to G_{2} \to H_{2} \to 1
\]
where $H_{1}$ is linear and $H_{2}$ is quadratic.
\end{cor}
{\bf Proof.}

We let $G_{3}(X)$ be the image of the natural map $H^{2}(X, \mathcal{O}^{\times}) \to H^{0}(S,R^{2}f_{*}\mathcal{O}^{\times})$. Let $G$ be the kernel of the map to $F \to G_{3}$, $G_{0}(X) = H^{2}(S, \mathcal{O}^{\times})$, and
$G_{2}(X) = H^{1}(S, R^{1} f_{*}\mathcal{O}^{\times})$.
Let $K_{2}(X)$ and $K_{3}(X)$ be the images of the natural maps
\[H^{0}(S, R^{2} f_{*}\mathcal{O}) \cap G_{3}(X) \to G_{3}(X) \to H^{0}(S, R^{3} f_{*} \mathbb{Z}).
\]
Let $H_{1}(X)$ and $H_{2}(X)$ be the images of the natural maps
\[H^{1}(S, R^{1} f_{*}\mathcal{O}) \to
H^{1}(S, R^{1} f_{*}\mathcal{O}^{\times})
\to H^{1}(S, R^{2} f_{*} \mathbb{Z}).
\]

\ \hfill $\Box$

Let $A$ be an complex analytic abelian group space over $S$.  We will
always consider its section to be the identity of the relative
group structure over $S$.
Let
\[A^{n} = A_{(0)} \times_{S} \cdots \times_{S}A_{(n-1)}
\]
denote the $n-$fold fiber product where each $A_{(j)} = A$.
We will say that $A$ satisfies the generalized theorem of the $n-$cube
with respect to a functor $F$ as in the above corollary if the map
\[F(\sigma_{i}^{n}):F( A^{n}) \to F(A^{n-1})
\]
is injective.  In other words $A$ satisfies the generalized theorem of
the $n-$cube with respect to $F$ precisely when the functor $F$ is of
order $n-1$ when restricted to the object $A$ and its
powers.  Following tradition, we call the generalized theorem
of the $1-$cube the generalized theorem of the segment, the generalized
theorem of the $2-$cube will be
called the generalized theorem of the square, the generalized theorem
of the $3-$cube will be called
the generalized theorem of the cube and the generalized theorem of
the $4-$cube will be called the generalized
theorem of the hyper-cube. Theorem \ref{thm:CubicFunctor} shows that complex tori $A \to S$ in the 
relative sense satisfy the generalized theorem
of the hyper-cube with respect to the functor $F(A) = H^{2}(A,\mathcal{O}^{\times})$.  This result parallels a result in \cite{Ho1972} and can be
used in conjunction with Corollary \ref{cor:filtration} to show that any torsion gerbe on $A$ is representable by an Azumaya Algebra \cite{Be1972}, \cite{ElNa1983}, \cite{Ho1972}.

Let us derive some consequences of Theorem \ref{thm:CubicFunctor}.   Let $A$ be an complex analytic abelian group space over $S$.  For any four maps $f_{i}:X_{i} \to A$ over $S$, such that $f_{i}|_{S}$ is the identity section of $A$ we see that there is an isomorphism of gerbes on $X_{1} \times_{S} X_{2} \times_{S} X_{3} \times_{S} X_{4}$
\begin{equation}\label{eqn:FourThings}
\begin{split}
& \Big{(}(f_{1}+f_{2}+f_{3}+f_{4})^{*} \mathfrak{G}\Big{)} \otimes
\Big{(}\prod_{1 \leq i<j \leq 4}
(f_{i}+f_{j})^{*}\mathfrak{G} \Big{)} \\
& \cong \Big{(} \prod_{1 \leq i<j <k \leq 4}
(f_{i}+f_{j}+f_{k})^{*}\mathfrak{G} \Big{)} \otimes \Big{(} \prod_{1 \leq i\leq 4}
f_{i}^{*}\mathfrak{G} \Big{)}
\end{split}
\end{equation}
for the pullbacks of a gerbe $\mathfrak{G}$ on $A$ under of the various maps \[X_{1} \times_{S} X_{2} \times_{S} X_{3} \times_{S} X_{4} \to A.\]
Indeed a simple calculation shows that the left and right side of (\ref{eqn:FourThings}) are isomorphic after restriction to any three out of the four spaces.  Let $a,b,c$ be sections of $A \to S$ and $X_{1}=A$, $X_{2}=X_{3}=X_{4}=S$, let $f_{1}$ be the identity map, and let $f_{2},f_{3},f_{4}$ be the sections $a$, $b$, and $c$ respectively.  Now as a consequence of (\ref{eqn:FourThings}), we have the following corollary of theorem \ref{thm:CubicFunctor}.
 \begin{cor}
 For any gerbe $\mathfrak{G}$ on an complex analytic abelian group space $A \to S$ and any sections $a,b,c$ of $A\to S$
\begin{equation}\label{eqn:TransFromMult}
t_{a+b+c}^{*}\mathfrak{G} \cong \mathfrak{G} \otimes
t_{a}^{*} \mathfrak{G}^{-1} \otimes
t_{b}^{*} \mathfrak{G}^{-1} \otimes
t^{*}_{c} \mathfrak{G}^{-1} \otimes
t_{a+b}^{*} \mathfrak{G} \otimes
t_{b+c}^{*} \mathfrak{G} \otimes
t^{*}_{a+c} \mathfrak{G}
\end{equation}
\end{cor}
Similarly, if $\mathfrak{G}$ is topologically trivial on the fibers, there are equations analogous to (\ref{eqn:FourThings}) and (\ref{eqn:TransFromMult}) involving three spaces $X_{i}$ or two sections of $A$ respectively.

Let $X$ be a complex torus and let $\mathfrak{G}$ be a gerbe on $X$. If $n$ is any integer, we denote by $n$ the corresponding isogeny $X \to X$.  We would like to say something about the pullback $n^{*} \mathfrak{G}$.  Our Appell-Humbert theorem shows that $\mathfrak{G}$ is a gerbe isomorphic to that given (\ref{eqn:ob2CyclesToGerbes}) by a cocycle $\Phi^{(B,E)} \in Z^{2}(\Lambda, \mathcal{O}^{\times}(V))$ described in definition \ref{dfn:Reps}.  An easy computation shows that
\[(-1)^{*} \Phi^{(B,E)} = \Phi^{(B,-E)}
\]
and
\[n^{*} \Phi^{(B,E)} = \Phi^{(n^{2}B,n^{3}E)}.
\]
and therefore we have the following corollary of remark \ref{rmk:important}.
\begin{cor}
For any integer $n$ and any gerbe $\mathfrak{G}$ on a complex torus $X$, there is an isomorphism
\[n^{*} \mathfrak{G} \cong \Big{(} \mathfrak{G}^{(\frac{n^{2} + n^{3}}{2})}\Big{)} \otimes \Big{(}(-1)^{*} \mathfrak{G}^{(\frac{n^{2} - n^{3}}{2})}\Big{)}
.\]
\end{cor}
This reduces to
\[n^{*} \mathfrak{G} \cong \mathfrak{G}^{(n^{2})}
\]
in the topologically trivial case.  The same formulas hold for $\mathfrak{G}$ replaced by an element in $Q(A,\varpi)$, see (\ref{eqn:DefOfQ}) for $\varpi: A \to S$ a complex analytic abelian group space.  There the topological triviality option is of course replaced by fiber-wise topological triviality.  One can use an induction argument to see that this is consistent with (\ref{eqn:FourThings}) in the case that all $f_{i}$ are replaced by isogenies.

\section{The Universal Gerbe}\label{univ}

In this section we find a complex analytic stack which serves as a fine moduli stack for topologically trivial holomorphic gerbes on a complex torus.  This is analogous to $\text{Pic}^{0}(X) = \overline{V}^{\vee}/\Lambda^{\vee}$ for topologicaly trivial line bundles.  The short exact sequence (\ref{exponentialSES})
suggests that this stack could be $[H^{2}(X, \mathcal{O})/H^{2}(X, \mathbb{Z}) ]$.  Accepting this suggestion, our task becomes to find a universal gerbe $\mathfrak{P}$ on $X \times [H^{2}(X, \mathcal{O})/H^{2}(X, \mathbb{Z}) ]$.  We conclude this section by calculating the topological type of this universal gerbe.  For some general background material on gerbes over stacks see \cite{FeHeRoZh2008}, \cite{Po2008}, \cite{BeXu2003} and the references in those papers.  In order to do this we parameterize this stack as $[\wedge^{2} \overline{V}^{\vee} / \text{Alt}^{2}(\Lambda, \mathbb{Z})]$ where $Alt^{2}(\Lambda, \mathbb{Z})$ acts by its image subgroup \[\text{Alt}^{2}(\Lambda, \mathbb{Z})^{H} \subset \wedge^{2} \overline{V}^{\vee}\] under the Hodge projection
\[H^{2}(X, \mathbb{Z}) \subset H^{2}(X,\mathbb{R}) \to H^{2}(X, \mathbb{C}) \to H^{2}(X, \mathcal{O})=H^{0,2}(X)
\]
or
\[
\text{Alt}^{2}(\Lambda, \mathbb{Z}) \subset \text{Alt}^{2}(\Lambda, \mathbb{R}) \to \text{Alt}^{2}(\Lambda, \mathbb{C}) \to \wedge^{2}\overline{V}^{\vee}.
\]
Explicitly, for $v_{1}, v_{2} \in V$, the Hodge projection map takes $\omega$ to $\omega^{H}$ where
\begin{equation}\label{eqn:ActualHodgeProj}
\omega^{H}(v_{1}, v_{2}) = \frac{1}{4}\bigg{(}\omega(v_{1}, v_{2})-\omega(iv_{1}, iv_{2})+i\omega(iv_{1}, v_{2})+i\omega(v_{1}, iv_{2})\bigg{)}.
\end{equation}
In contrast to the situation in section \ref{lineAH},  the image of $\text{Alt}^{2}(\Lambda, \mathbb{Z})$ might not be closed inside $\wedge^{2}\overline{V}^{\vee}$ and cannot be described by equations.

We will describe a gerbe on the stack
\[V/\Lambda \times [\wedge^{2} \overline{V}^{\vee} / \text{Alt}^{2}(\Lambda,\mathbb{Z})] = [(V \times \wedge^{2} \overline{V}^{\vee})/( \Lambda \times \text{Alt}^{2}(\Lambda,\mathbb{Z}) )]
\]
by writing down an element
\[\Psi \in Z^{2}(\Lambda \times \text{Alt}^{2}(\Lambda, \mathbb{Z}), \mathcal{O}^{\times}(V \times \wedge^{2}\overline{V}^{\vee}))\]  and associating to it a gerbe via the method explained in the second appendix.  In fact all gerbes on this stack come about from this method, see  Remark \ref{rem:GerbeOnStackFromGroupCoho} and the second appendix in Section \ref{appendix2}.  We will then show that restricting the gerbe to a point of the stack, gives a gerbe whose associated class in $H^{2}(\Lambda, \mathcal{O}^{\times}(V))$ equals the cohomology class of the image of the point under $\sigma$.

We define below a map
\[\Psi: (\Lambda \times \text{Alt}^{2}(\Lambda, \mathbb{Z}))^{2} \to \mathcal{O}^{\times}(V \times \wedge^{2}\overline{V}^{\vee}).
\]
Given elements $v \in V$, $B \in \wedge^{2}\overline{V}^{\vee}$, $\lambda_{1}, \lambda_{2} \in \Lambda$ and $\mu_{1}, \mu_{2} \in \text{Alt}^{2}(\Lambda, \mathbb{Z})$ let
\begin{equation}\label{eqn:DefOfPoincarePsi}
\begin{split}
\Psi_{(\lambda_{1},\mu_{1}),(\lambda_{2}, \mu_{2})}(v,B) = \exp \bigg{(} & \frac{1}{2}(B+\mu_{1}^{H} +\mu_{2}^{H} )(\lambda_{1}, \lambda_{2})+\frac{1}{2}\overline{\mu_{2}^{H}(v, \lambda_{1})} \\
& -\frac{i}{4}\Big{(}\mu_{2}(iv,\lambda_{1})-\mu_{2}(v,i\lambda_{1})\Big{)} \\
& + \frac{1}{4}\Big{(}\mu_{2}(v,\lambda_{1}) +\mu_{2}(iv,i\lambda_{1})\Big{)} \\
&-\frac{i}{8}\Big{(}\mu_{2}(i\lambda_{1}, \lambda_{1})-\mu_{2}(\lambda_{1},i\lambda_{1})\Big{)} \\
&+\frac{1}{2}\sigma(\mu_{2})(\lambda_{1},\lambda_{1})\bigg{)} .
\end{split}
\end{equation}

The map $\sigma$ is defined in (\ref{eqn:DefOfSigma}). It is easily seen that $\Psi_{(\lambda_{1},\mu_{1}),(\lambda_{2}, \mu_{2})}$ is a holomorphic
function on $V \times \wedge^{2} \overline{V}^{\vee}$
for every value of
$\lambda_{1}, \lambda_{2}, \mu_{1}, \mu_{2}$.  Lest
this equation look too mysterious, let the
reader notice that the terms on the second, third,
and fourth line have same form as part of the
canonical factor of automorphy for line bundles
where the form of type $(1,1)$ has become
\[\frac{1}{2}\Big{(}\mu_{2}(v_{1},v_{2})
+\mu_{2}(iv_{1},iv_{2})\Big{)},
\]
although this is not integer valued on
the lattice.  We need to show that the boundary

\begin{equation}
\begin{split}
& (\delta \Psi)_{(\lambda_{1}, \mu_{1}), (\lambda_{2}, \mu_{2}), (\lambda_{3}, \mu_{3})}(v,B) \\
& = (\Psi_{(\lambda_{2}, \mu_{2}), (\lambda_{3}, \mu_{3})}(v+\lambda_{1},B+\mu_{1}^{H}))  (\Psi_{(\lambda_{1}+\lambda_{2}, \mu_{1}+\mu_{2}), (\lambda_{3}, \mu_{3})}(v,B))^{-1} \\ & \quad (\Psi_{(\lambda_{1}, \mu_{1}),(\lambda_{2}+\lambda_{3}, \mu_{2}+ \mu_{3})}(v,B))  (\Psi_{(\lambda_{1}, \mu_{1}),(\lambda_{2}, \mu_{2})}(v,B))^{-1}
\end{split}
\end{equation}
is trivial.  The terms on the first line of (\ref{eqn:DefOfPoincarePsi}) have boundary
\[\exp \left(  \frac{1}{2}\mu_{3}^{H}(\lambda_{1}, \lambda_{2})+\frac{1}{2}\overline{\mu_{3}^{H}(\lambda_{1}, \lambda_{2})} \right)
\]
The terms on the second, third, and fourth line of (\ref{eqn:DefOfPoincarePsi}) have boundary
\[\exp \left(\frac{1}{4}\Big{(}\mu_{3}(\lambda_{1},\lambda_{2})
+\mu_{3}(i\lambda_{1},i\lambda_{2})\Big{)}\right).
\]
  The boundary term on the last line of (\ref{eqn:DefOfPoincarePsi}) is
\begin{equation}
\begin{split}
& \exp \left(-\frac{1}{2}\Big{(}\sigma(\mu_{3})(\lambda_{1},\lambda_{2})
+\sigma(\mu_{3})(\lambda_{2},\lambda_{1})\Big{)}\right) \\
& = \exp\left(-\frac{1}{2}(s \sigma (\mu_{3}))(\lambda_{1},\lambda_{2})\right) = \exp\left(\frac{1}{2}\mu_{3}(\lambda_{1},\lambda_{2}) \right)
\end{split}
\end{equation}
where we have used $\exp(\frac{1}{2}\mathbb{Z}) \subset \{\pm 1\}$.
Therefore using (\ref{eqn:ActualHodgeProj})
\begin{equation}
\begin{split}
& (\delta \Psi)_{(\lambda_{1}, \mu_{1}), (\lambda_{2}, \mu_{2}), (\lambda_{3}, \mu_{3})}(v,B) \\
&  =
\exp \left(  \frac{1}{2}\mu_{3}^{H}(\lambda_{1}, \lambda_{2})+\frac{1}{2}\overline{\mu_{3}^{H}(\lambda_{1}, \lambda_{2})}   +\frac{1}{4}\Big{(}\mu_{3}(\lambda_{1},\lambda_{2})
+\mu_{3}(i\lambda_{1},i\lambda_{2})\Big{)} +\frac{1}{2}\mu_{3}(\lambda_{1}, \lambda_{2}) \right) \\
& = \exp\left( \frac{1}{4}\Big{(}\mu_{3}(\lambda_{1},\lambda_{2})-\mu_{3}(i\lambda_{1},i\lambda_{2})\Big{)}
+\frac{1}{4}\Big{(}\mu_{3}(\lambda_{1},\lambda_{2})+\mu_{3}(i\lambda_{1},i\lambda_{2})\Big{)}+\frac{1}{2}\mu_{3}(\lambda_{1}, \lambda_{2})\right) \\
& =  \exp \left( \mu_{3}(\lambda_{1},\lambda_{2}) \right) \\
& = 1.
\end{split}
\end{equation}
\begin{dfn}
Let $\mathfrak{P}$ be the gerbe on $X \times [H^{2}(X,\mathcal{O})/H^{2}(X,\mathbb{Z})]$ corresponding under (\ref{eqn:ob2CyclesToGerbes}) to $\Psi$.
\end{dfn}
Restricting the gerbe $\mathfrak{P}$ to $X \times [B ]$ gives a gerbe on $X$ isomorphic to that given by the cocycle in $Z^{2}(\Lambda, \mathcal{O}^{\times}(V))$
\[\exp \left( \frac{1}{2} B(\lambda_{1}, \lambda_{2}) \right)
\]
whose equivalence class is the exponential of
\[B \in \wedge^{2}\overline{V}^{\vee} \cong H^{2}(X,\mathcal{O}) \to H^{2}(X,\mathcal{O}^{\times}).
\]
\begin{thm}\label{thm:GerbeFine}
The gerbe $\mathfrak{P}$ on $X \times [H^{2}(X \times \mathcal{O})/H^{2}(X,\mathbb{Z})]$ has the property that for any connected complex analytic space $T$ and any gerbe $\mathfrak{G} \to X \times T$ which is topologically trivial on each fiber $X_{t}$, there is a holomorphic map
\[f: T \to [H^{2}(X , \mathcal{O})/H^{2}(X,\mathbb{Z})]
\] such that
\[\mathfrak{G} \cong ((1,f)^{*}\mathfrak{P}) \otimes \mathfrak{C}\] where $\mathfrak{C}$ is a gerbe trivial on each fiber.
\end{thm}
\noindent
Notice that we have already proven this in the case that $T$ is a point.


{\bf Proof.}
Let $U_{T}$ be the universal cover of $T$.  We will denote by \[\rho:X \times T \to T,\] and \[\tilde{\rho}: X \times U_{T} \to U_{T},\] and \[ \tilde{\tilde{\rho}}: X \times \pi_{1}(T) \times U_{T} \to  \pi_{1}(T) \times U_{T} ,\] the natural projections.  Let 
\[G \in H^{0}(U_{T},(R^{2} \tilde{\rho}_{*} \mathcal{O}^{\times})_{0})
\]
be defined by pulling back $\mathfrak{G}$ to $X \times U_{T}$ and then restricting it to the fibers of the projection to $U_{T}$.
Define
\[(R^{2} \rho_{*} \mathcal{O}^{\times})_{0} = Im[R^{2} \rho_{*}\mathcal{O} \to R^{2} \rho_{*}\mathcal{O}^{\times}] = \text{ker}[R^{2} \rho_{*}\mathcal{O}^{\times} \to R^{3} \rho_{*} \mathbb{Z}]
\]
and
\[(R^{2} \tilde{\rho}_{*} \mathcal{O}^{\times})_{0} = Im[R^{2} \tilde{\rho}_{*}\mathcal{O} \to R^{2} \tilde{\rho}_{*}\mathcal{O}^{\times}] = \text{ker}[R^{2} \tilde{\rho}_{*}\mathcal{O}^{\times} \to R^{3} \tilde{\rho}_{*} \mathbb{Z}]
.\]

We now consider the short exact sequence
\[1 \to TL_{X}
\to R^{2} \tilde{\rho}_{*}\mathcal{O} \to (R^{2} \tilde{\rho}_{*} \mathcal{O}^{\times})_{0} \to 1.
\]
of sheaves on $U_{T}$, where $TL_{X}$ is the transcendental lattice of $X$ thought of as a constant sheaf on $U_{T}$.
The obstruction to lifting our element $G \in H^{0}(U_{T},(R^{2} \tilde{\rho}_{*} \mathcal{O}^{\times})_{0})$ to $H^{0}(U_T,R^{2} \tilde{\rho}_{*}\mathcal{O})$ lives in $H^{1}(U_{T}, TL_{X})=0$.  And so such a lift is always possible.  The term $H^{0}(U_{T},R^{2}\tilde{\rho}_{*} \mathcal{O})$ represents the holomorphic maps from $U_{T}$ to $\wedge^{2}\overline{V}^{\vee}$
\begin{equation} \label{eqn:rep1}
H^{0}(U_{T},R^{2}\tilde{\rho}_{*} \mathcal{O}) = Hol(U_{T}, \wedge^{2} \overline{V}^{\vee})
.
\end{equation}
Define \[f_{\mathfrak{G}, \mathbb{C}}:U_{T} \to \wedge^{2} \overline{V}^{\vee}\] to be any lift of $G$.


We can calculate the group cohomology of $\pi_{1}(T)$ acting on the global sections of any sheaf of groups $\mathcal{S}$ on $U_{T}$ by taking the homology groups of the sequence
\[H^{0}(U_{T},\mathcal{S}) \to H^{0}(\pi_{1}(T) \times U_{T}, \mathcal{S}) \to  H^{0}(\pi_{1}(T) \times \pi_{1}(T) \times U_{T}, \mathcal{S}) \to \cdots.
\]
We can combine these with the long exact sequences induced by the short exact sequence of sheaves on $U_{T}$
\[0 \to TL_{X} \to R^{2} \tilde{\rho}_{*} \mathcal{O} \to (R^{2} \tilde{\rho}_{*} \mathcal{O}^{\times})_{0} \to 1.
\]
The result is a commutative diagram with exact columns
\begin{equation}\label{eqn:DiagCompat}
\entrymodifiers={+!!<0pt,\fontdimen22\textfont2>}
\xymatrix @W-15pc @C-.4pc {\ar @{} [dr] |{}
&H^{0}(T, (R^{2} \rho_{*} \mathcal{O}^{\times})_0) \ar@{^{(}->}[r]
& H^{0}(U_T,(R^{2}\tilde{\rho}_{*}\mathcal{O}^{\times})_0)
\ar[r] &
H^{0}(\pi_{1}(T) \times U_T, (R^{2}\tilde{\tilde{\rho}}_{*} \mathcal{O}^{\times})_0) \\
& H^{0}( T, \wedge^{2}\overline{V}^{\vee} \otimes \mathcal{O}) \ar[u] \ar@{^{(}->}[r]  &
H^{0}(U_T,\wedge^{2}\overline{V}^{\vee} \otimes \mathcal{O})
\ar@{>>}[u] \ar[r] &
H^{0}(\pi_{1}(T) \times U_T,\wedge^{2}\overline{V}^{\vee} \otimes \mathcal{O})  \ar[u]  \\
  &
 H^{0}(T,TL_{X}) \ar@{^{(}->}[u] \ar [r]^{\cong}  & H^{0}(U_T , TL_{X}) \ar@{^{(}->}[u]  \ar[r]^{0} &
 H^{0}(\pi_{1}(T) \times U_{T},TL_{X}) \ar[r] \ar@{^{(}->}[u]   &}
\end{equation}
The surjectivity is a result of $H^{1}(U_T, TL_{X}) = 0$.  A simple diagram chase establishes the a {\it diagonal map}
\[D: H^{0}(T, (R^{2}\rho_{*} \mathcal{O}^{\times})_0) \to H^{0}(\pi_{1}(T) \times U_{T},TL_{X}).\]  In terms of the diagram (\ref{eqn:DiagCompat}), it corresponds to moving right, down, right, down.  In fact, the image of this map is in the kernel of the map
\[H^{0}(\pi_{1}(T) \times U_{T},TL_{X})  \to H^{0}(\pi_{1}(T) \times \pi_{1}(T) \times U_{T},TL_{X}).
\]  Upon taking homology in the horizontal direction the snake like map corresponds to the first boundary map in group cohomology
\begin{equation} \label{eqn:GivesInt}
\xymatrix  @W-2pc @C-1.1pc
{H^{0}(T, (R^{2}\rho_{*} \mathcal{O}^{\times})_0) =
H^{0}(U_T,(R^{2}\tilde{\rho}_{*} \mathcal{O}^{\times})_0)^{\pi_{1}(T)} \ar[r]
& H^{1}(\pi_{1}(T), TL_{X}) = Hom(\pi_{1}(T), TL_{X})}\end{equation}
 coming from the short exact sequence
\[
\xymatrix   { 0 \ar[r] & TL_{X} \ar[r] & H^{0}(U_T, \wedge^{2} \overline{V}^{\vee} \otimes \mathcal{O}) \to H^{0}(U_T,  (R^{2}\tilde{\rho}_{*} \mathcal{O}^{\times})_0) \ar[r]  & 1
}
\]
of $\pi_{1}(T)$ modules.
The diagonal map (\ref{eqn:GivesInt}) factors using the Hodge projection
\begin{equation}\label{eqn:something_else}
\text{Alt}^{2}(\Lambda, \mathbb{Z}) = H^{2}(X,\mathbb{Z}) \to TL_{X}\end{equation}
\[\mu \mapsto \mu^{H}
\] as
\[\xymatrix{H^{0}(T, (R^{2}\rho_{*} \mathcal{O}^{\times})_0) \ar[r] & Hom(\pi_{1}(T), \text{Alt}^{2}(\Lambda, \mathbb{Z})) \ar[r]  & Hom(\pi_{1}(T), TL_{X}).}
\]
In terms of the diagram (\ref{eqn:DiagCompat}) we took $f_{\mathfrak{G}, \mathbb{C}}$ to be any element of $H^{0}(U_{T}, \wedge^{2} \overline{V}^{\vee} \otimes \mathcal{O})$ covering $G$.  Consider the
image 
of $G$ under the map on cohomology $H(D)$ induced by the diagonal map $D$ into
\[Hom(\pi_{1}(T), TL_{X}) = \text{ker}\Big{[}H^{0}(\pi_{1}(T) \times U_{T},TL_{X}) \to H^{0}(\pi_{1}(T) \times \pi_{1}(T) \times U_{T},TL_{X})\Big{]}.
\]
By choosing a group theoretic splitting of (\ref{eqn:something_else}) we can chose lifts of such elements to 
\[f_{\mathfrak{G},\mathbb{Z}} \in \text{Hom}(\pi_{1}(T), \text{Alt}^{2}(\Lambda, \mathbb{Z})).
\]
Therefore 
\[H(D)(G) = f_{\mathfrak{G},\mathbb{Z}} ^{H}
.\]
Finally we conclude \begin{equation}\label{eqn:compat}f_{\mathfrak{G}, \mathbb{C}}(\chi \cdot u) = (f_{\mathfrak{G}, \mathbb{Z}}(\chi))^{H} + f_{\mathfrak{G}, \mathbb{C}}(u)
\end{equation}
for all $\chi \in \pi_{1}(T)$ and $u \in U_{T}$.  This  follows from the fact that $f_{\mathfrak{G}, \mathbb{C}}$ and
 $f_{\mathfrak{G}, \mathbb{Z}}^{H}$ map to the same middle term of the right-most
column of (\ref{eqn:DiagCompat}).

Together, $f_{\mathfrak{G}, \mathbb{Z}}$ and $f_{\mathfrak{G}, \mathbb{C}}$ define a holomorphic map
\[f: T \cong U_{T} / \pi_{1}(T) \to [\wedge^{2}\overline{V}^{\vee} / \text{Alt}^{2}(\Lambda, \mathbb{Z}) ].
\]
In what may perhaps be more familiar terms to the reader, we can rearrange this data into an equivalent description given by a diagram
\[
\xymatrix{\ar @{} [dr] |{}
& (U_{T}) \times_{\pi_{1}(T)} \text{Alt}^{2}(\Lambda, \mathbb{Z}) \ar[d] \ar[r]^-{\Upsilon} & \wedge^{2}\overline{V}^{\vee} &  \\
  &
 T  &   & .}
\]
In this diagram the vertical arrow is an $\text{Alt}^{2}(\Lambda, \mathbb{Z})$ principal bundle on $T$ induced via $f_{\mathfrak{G}, \mathbb{Z}}$ by the $\pi_{1}(T)$ principal bundle $U_{T} \to T$.   The horizontal arrow is an $\text{Alt}^{2}(\Lambda, \mathbb{Z})$-equivariant map given by to adding the map $f_{\mathfrak{G}, \mathbb{C}}$ to the Hodge projection map.
\[\Upsilon(u,\mu) = f_{\mathfrak{G}, \mathbb{C}}(u) + \mu^{H}
\]
This is well defined because of equation (\ref{eqn:compat}).
The gerbe $(1,f)^{*}\mathfrak{P}$ has cocycle representative \[(1,f)^{*}\Psi \in Z^{2}(\Lambda \times \pi_{1}(T),\mathcal{O}^{\times}(V \times U_{T})))\] given by
\[((\lambda_{1}, \chi_{1}),(\lambda_{2},\chi_{2})) \mapsto \Psi_{((\lambda_{1},f_{\mathfrak{G},\mathbb{Z}}(\chi_{1})),( \lambda_{2},  f_{\mathfrak{G},\mathbb{Z}}(\chi_{2}))}
(v,f_{\mathfrak{G},\mathbb{C}}(u))
\]
for $\lambda_{1}, \lambda_{2} \in \Lambda$, $\chi_{1}, \chi_{2} \in \pi_{1}(T)$ , $u \in U_{T}$ and $v \in V$.  Restricting $(1,f)^{*}\mathfrak{P}$ to $X \times \{t \}$ clearly gives a gerbe isomorphic $\mathfrak{G} |_{X \times \{ t \}}$. Let
\[\mathfrak{C} = ((1,f)^{*}\mathfrak{P})^{-1} \mathfrak{G}
\] be the correction term.  The restriction of $\mathfrak{G}$ to a point $t \in T$ gives a gerbe isomorphic to the gerbe corresponding to the cocycle
\[(\lambda_{1},\lambda_{2}) \mapsto \exp\left(\frac{1}{2}f_{\mathfrak{G},\mathbb{C}}(u)(\lambda_{1},\lambda_{2})\right).
\]
where $u$ is any lift of $t$. 
\ \hfill $\Box$

Let $G_{X}$ denote the action groupoid

\[
\xymatrix{**[l]{H^{2}(X,\mathbb{Z}) \times H^{2}(X, \mathcal{O})}  \ar@/_1.5pc/[r]_-{
} \ar@/^1.5pc/[r]^-{
} & **[r] \ar@/^0pc/[l]_-{}  H^{2}(X, \mathcal{O}) }
\]
Notice that the above theorem provided a map from the objects of $Hol(T,G_{X})$ to gerbes on $X \times T$ which are topologically trivial on each fiber.  We can in fact promote this to a functor.  The functor is a combination of the previous theorem together with the following definition which defines the functor on morphisms.  Because $T$ is connected and the image of $H^{2}(X,\mathbb{Z})$ in $H^{2}(X, \mathcal{O})$ is countable, there is a morphism (isomorphism) between the gerbes $(1,f)^{*} \mathfrak{P}$ and $(1,g)^{*} \mathfrak{P}$ if and only if the difference $f - g$ comes from some $\mu \in \text{Alt}^{2}(\Lambda, \mathbb{Z})$.  The space of morphisms between $f$ and $g$ is the group of $\mu$ such that $\mu^{H} = f-g$.  Hence we need to give a functorial assignment to any such $\mu$ of a trivialization of the gerbe $(1,\mu^{H})^{*}\mathfrak{P}$.  This gerbe by defintion is the gerbe coming from the cocycle
\[(\lambda_{1}, \lambda_{2}) \mapsto \exp \left(\frac{1}{2}\mu^{H}(\lambda_{1}, \lambda_{2}) \right)
.\]
\begin{lem}
For every $\mu \in \textup{Alt}^{2}(\Lambda, \mathbb{Z})$ we can assign a global trivialization $\mathcal{L}_{\mu}$ of $(1, \mu^{H})^{*} \mathfrak{P}$ in a way that
\[\mathcal{L}_{\mu} \otimes \mathcal{L}_{\mu'} \cong \mathcal{L}_{\mu+\mu'}.
\]
and these isomorphisms are compatible in the obvious way for three values of $\textup{Alt}^{2}(\Lambda, \mathbb{Z})$.
\end{lem}
{\bf Proof.}
We need to define, in a way linear in $\mu$, a cochain in $C^{1}(\Lambda \times \pi_{1}(T), \mathcal{O}^{\times}(V \times U_{T}))$ whose boundary in $C^{2}(\Lambda \times \pi_{1}(T), \mathcal{O}^{\times}(V \times U_{T}))$ is
\[(\lambda_{1}, \lambda_{2}) \mapsto \exp \left(\frac{1}{2}\mu^{H}(\lambda_{1}, \lambda_{2}) \right).
\]
The type decomposition presents any $\mu$ as sum of alternating pieces
\[\mu = \mu^{(0,2)} + \mu^{(1,1)} + \mu^{(2,0)} = \mu^{H} + \mu^{(1,1)} + \overline{\mu^{H}}.
\]
It is easy to bound separately the terms aside from $\mu^{H}$ and combine them into a boundary for
\[\frac{1}{2}\mu^{H}  = \frac{1}{2}\Big{(}[\mu] -  [\mu^{(1,1)}] - [\overline{\mu^{H}}]\Big{)}.\]  An easy calculation shows that the map
\[\kappa:\Lambda \times \pi_{1}(T) \to \mathcal{O}^{\times}(V \times U_{T})
\]
defined by
\begin{equation}\label{eqn:IsomInGerbes}
\begin{split}
& (\lambda, \chi) \mapsto \\ & \exp \left(\frac{1}{2}\Big{(}\Big{[}\sigma(\mu)(\lambda, \lambda)\Big{]}-\Big{[}\frac{-i}{2}\mu^{(1,1)}(iv,\lambda)+\frac{1}{2}\mu^{(1,1)}(v,\lambda)
-\frac{i}{4}\mu^{(1,1)}(i\lambda,\lambda)\Big{]}
-\Big{[}\overline{\mu^{H}(v,\lambda)}\Big{]}\Big{)}\right)
\end{split}
\end{equation}
does the job.  We define $\mathcal{L}_{\mu}$ to be the trivialization of $(1,\mu^{H})^{*} \mathfrak{P}$ corresponding (\ref{eqn:2CyclesToGerbes}) to $\kappa$ as explained in the second appendix in section \ref{appendix2}.
\ \hfill $\Box$
\begin{cor}
Let $\mathfrak{G}erbes(X \times T)'$ be the category of gerbes topologically trivial on each fiber of the projection to $T$.  We have defined a functor
\[Hol(T, [H^{2}(X,\mathcal{O})/H^{2}(X,\mathbb{Z})]) \to \mathfrak{G}erbes(X \times T)'.
\]
The corresponding map on equivalence classes
\[\frac{Hol\Big{(}T, [H^{2}(X,\mathcal{O})/H^{2}(X,\mathbb{Z})]\Big{)}}{\sim} \quad \quad  \cong \quad \quad H^{0}\Big{(}T, (R^{2} \rho_{*} \mathcal{O}^{\times})_0\Big{)}
\]
is an isomorphism of groups.  Thus every gerbe on $X \times T$ topologically trivial along the fibers is equivalent up to multiplication by gerbes trivial along the fibers to a pullback of $\mathfrak{P}$ under a map $(1 \times f)$ where $f$ is unique up to equivalence.  This, in combination with Theorem \ref{thm:GerbeFine} concludes the proof that
$[H^{2}(X, \mathcal{O})/H^{2}(X, \mathbb{Z})]$ is a fine moduli stack.
\end{cor}

\ \hfill $\Box$


We now compute the topological type of $\mathfrak{P}$.  The topological class of $\mathfrak{P}$ lives in
\[H^{3}\Big{(}X \times [H^{2}(X , \mathcal{O})/H^{2}(X,\mathbb{Z})], \mathbb{Z}\Big{)} = H^{3}\Big{(}\Lambda \times \text{Alt}^{2}(\Lambda, \mathbb{Z}), \mathbb{Z}\Big{)}.
\]
The cocycle $\delta \log \Psi \in Z^{3}(\Lambda \times \text{Alt}^{2}(\Lambda, \mathbb{Z}),\mathbb{Z})$ is
\[(\delta \log \Psi)_{(\lambda_{1}, \mu_{1}),(\lambda_{2},\mu_{2}),(\lambda_{3}, \mu_{3})} =
\frac{1}{2}\mu_{3}(\lambda_{1}, \lambda_{2}) -\frac{1}{2}\Big{(}\sigma(\mu_{3})(\lambda_{1},\lambda_{2})
+\sigma(\mu_{3})(\lambda_{2},\lambda_{1})\Big{)}.\]
The associated class
\[s(\delta \log \Psi) \in \text{Alt}^{3}\Big{(}\Lambda \times \text{Alt}^{2}(\Lambda, \mathbb{Z}),\mathbb{Z}\Big{)}
\]
given as the image of the isomorphism
\[s:H^{3}\Big{(}\Lambda \times \text{Alt}^{2}(\Lambda, \mathbb{Z}), \mathbb{Z}\Big{)} \cong \text{Alt}^{3}\
\Big{(}\Lambda \times \text{Alt}^{2}(\Lambda, \mathbb{Z}),\mathbb{Z}\Big{)}
\]
is the skew-symmetrization of the above cocycle, which is easily calculated to be the assignment
\[(\lambda_{1}, \lambda_{2}, \lambda_{3}, \mu_{1}, \mu_{2}, \mu_{3}) \mapsto \mu_{1}(\lambda_{2}, \lambda_{3})-\mu_{2}(\lambda_{1}, \lambda_{3})+\mu_{3}(\lambda_{1}, \lambda_{2}),
\]
because the skew-symmetrization of $\sigma(\mu_{3})(\lambda_{1},\lambda_{2})+\sigma(\mu_{3})(\lambda_{2},\lambda_{1})$ is zero.
\section{Future Directions}
Perhaps the strongest missing link with the analogy we have drawn between line bundles and gerbes on complex tori is the question of theta functions.  What is the correct analogue of theta functions for the gerbes we have described?  There does not seem to exist a reasonable definition of {\em positivity} in this context, and hence allowing "zeros", no holomorphic sections exist.  One option is to study meromorphic sections of a given gerbe.  This direction has been pursued in a different context in the paper of Felder, Henriques, Rossi, and Zhu \cite{FeHeRoZh2008} and we hope that in our context those kind of relationships exist as well.  There are likely to be other connections with that work, for instance on a complex surface with Picard number $3$, our moduli stack is a trivial gerbe over a 'triplic curve' defined in \cite{FeHeRoZh2008} .  In future work, the author will explain various other phenomena including the representation by Azumaya algebras, and analogues of the theta and Heisenberg group.
\section{Appendix 1 }\label{appendix1}
Define the group $S(\Lambda)$ as the group of elements of $C^{2}(\Lambda, \mathbb{R}) \times A(\Lambda)$ which satisfy equation (\ref{eqn:real_part}).  We provide here a group theoretic splitting of the short exact sequence
\begin{equation}\label{eqn:GerbeSemiCharSplit}
1 \to Z^{2}(\Lambda, U(1)) \to S(\Lambda) \to A(\Lambda) \to 0.
\end{equation}
Fix an element $E \in A(\Lambda)$.  We wish to lift it to $(\beta',E) \in S(\Lambda)$.  Indeed the main purpose of this appendix is to find $\beta'$, the real part of $\beta$, which appears in the the cocycle 
\[\Phi^{E}_{\lambda_{1}, \lambda_{2}} = \exp(H_{\lambda_{1}, \lambda_{2}}(v) + \beta_{\lambda_{1}, \lambda_{2}})
\]
living in $Z^{2}(\Lambda, \mathcal{O}^{\times}(V))$ 
as in Definition \ref{DefOfGerbeForE}.  Notice that $\frac{1}{6}E(\lambda_{1},\lambda_{2}, \lambda_{3})$ agrees (\ref{eqn:skewK}) after skew-symmetrization with $k(\lambda_{1},\lambda_{2}, \lambda_{3})$ (defined in equation (\ref{eqn:keqn}))and hence by Lemma \ref{lem:skew} they are equivalent elements of $Z^{3}(\Lambda, \mathbb{R})$.  In order to find $\beta'$ we will first write down an element $u \in Map(\Lambda \times \Lambda, \mathbb{R})$ such that
\[
(\delta u)_{\lambda_{1},\lambda_{2}, \lambda_{3}} + \frac{1}{6}E(\lambda_{1},\lambda_{2}, \lambda_{3})= u_{\lambda_{2},\lambda_{3}}-u_{\lambda_{1}+\lambda_{2},\lambda_{3}}+u_{\lambda_{1},\lambda_{2}+\lambda_{3}}
-u_{\lambda_{1},\lambda_{2}}+\frac{1}{6}E(\lambda_{1},\lambda_{2}, \lambda_{3}) \in \mathbb{Z}.
\]
and then we will account for the difference between $\frac{1}{6}E(\lambda_{1},\lambda_{2}, \lambda_{3})$ and $k(\lambda_{1},\lambda_{2}, \lambda_{3})$. We can find $u$ by choosing an ordered basis $\{ \lambda^{i}\} $ for
$\Lambda$ and expanding the elements $\lambda_{\alpha} \in \Lambda$ as
$\lambda_{\alpha} = \sum_{i} n_{\alpha, i} \lambda^{i}$ for $n_{\alpha, i} \in \mathbb{Z}$.
\begin{equation}
\begin{split}
u_{\lambda_{1},\lambda_{2}} &= \sum_{i<j<k}\Big{[}\frac{5}{12}
E(n_{1,i}\lambda^{i},
n_{1,j}\lambda^{j},
n_{2,k}\lambda^{k})
+\frac{9}{12}
E(n_{1,i}\lambda^{i},
n_{2,j}\lambda^{j},
n_{1,k}\lambda^{k})
+\frac{1}{12}
E(n_{2,i}\lambda^{i},
n_{1,j}\lambda^{j},
n_{1,k}\lambda^{k}) \\
&  \quad +
\frac{5}{12}
E(n_{1,i}\lambda^{i},
n_{2,j}\lambda^{j},
n_{2,k}\lambda^{k})
+\frac{9}{12}
E(n_{2,i}\lambda^{i},
n_{1,j}\lambda^{j},
n_{2,k}\lambda^{k})
+\frac{1}{12}
E(n_{2,i}\lambda^{i},
n_{2,j}\lambda^{j},
n_{1,k}\lambda^{k})
\Big{]}.
\end{split}
\end{equation}
Using equation (\ref{eqn:helper}) we see that
\begin{equation}
\begin{split}
(\delta u)_{\lambda_{1},\lambda_{2}, \lambda_{3}}
& = \frac{5}{12}\Big{(}\sum_{i<j<k}
E(n_{1,i}\lambda^{i},
n_{2,j}\lambda^{j},
n_{3,k}\lambda^{k}) -
\sum_{j<i<k}E(n_{1,i}\lambda^{i},
n_{2,j}\lambda^{j},
n_{3,k}\lambda^{k})\Big{)} \\
 & \quad +\frac{9}{12}\Big{(}\sum_{i<k<j}
-E(n_{1,i}\lambda^{i},
n_{2,j}\lambda^{j},
n_{3,k}\lambda^{k})
+\sum_{j<k<i}E(n_{1,i}\lambda^{i},
n_{2,j}\lambda^{j},
n_{3,k}\lambda^{k})\Big{)} \\
 & \quad +\frac{1}{12}\Big{(}\sum_{k<i<j}
E(n_{1,i}\lambda^{i},
n_{2,j}\lambda^{j},
n_{3,k}\lambda^{k})
- \sum_{k<j<i}
E(n_{1,i}\lambda^{i},
n_{2,j}\lambda^{j},
n_{3,k}\lambda^{k})\Big{)} \\
 & \quad +\frac{5}{12}\Big{(}\sum_{i<j<k}
E(n_{1,i}\lambda^{i},
n_{2,j}\lambda^{j},
n_{3,k}\lambda^{k}) - \sum_{i<k<j}
E(n_{1,i}\lambda^{i},
n_{2,j}\lambda^{j},
n_{3,k}\lambda^{k})\Big{)}  \\
 & \quad +\frac{9}{12}\Big{(}\sum_{j<i<k}
-E(n_{1,i}\lambda^{i},
n_{2,j}\lambda^{j},
n_{3,k}\lambda^{k})
+ \sum_{k<i<j}E(n_{1,i}\lambda^{i},
n_{2,j}\lambda^{j},n_{3,k}\lambda^{k})\Big{)} \\
& \quad +\frac{1}{12}\Big{(}\sum_{j<k<i}
E(n_{1,i}\lambda^{i},n_{2,j}\lambda^{j},n_{3,k}\lambda^{k})
- \sum_{k<j<i}E(n_{1,i}\lambda^{i},
n_{2,j}\lambda^{j},n_{3,k}\lambda^{k})\Big{)}  \\
  &  = \frac{10}{12}\sum_{i<j<k}
  E(n_{1,i}\lambda^{i},
  n_{2,j}\lambda^{j},
  n_{3,k}\lambda^{k})
- \frac{14}{12}\sum_{i<k<j}
E(n_{1,i}\lambda^{i},n_{2,j}\lambda^{j},n_{3,k}\lambda^{k}) \\
& \quad -\frac{14}{12}\sum_{j<i<k}
E(n_{1,i}\lambda^{i},n_{2,j}\lambda^{j},n_{3,k}\lambda^{k})
+ \frac{10}{12}\sum_{j<k<i}
E(n_{1,i}\lambda^{i},n_{2,j}\lambda^{j},n_{3,k}\lambda^{k}) \\
 & \quad +\frac{10}{12}\sum_{k<i<j}
 E(n_{1,i}\lambda^{i},n_{2,j}\lambda^{j},n_{3,k}\lambda^{k})
-\frac{2}{12} \sum_{k<j<i}
E(n_{1,i}\lambda^{i},n_{2,j}\lambda^{j},n_{3,k}\lambda^{k}) \\
 & \equiv  -\frac{1}{6}\sum_{i,j,k}
 E(n_{1,i}\lambda^{i},n_{2,j}\lambda^{j},n_{3,k}\lambda^{k}) \\
& =  -\frac{1}{6}E(\lambda_{1},\lambda_{2}, \lambda_{3}).
\end{split}
\end{equation}
Under the decomposition
\[\beta' = u + r,
\]
all that remains is to solve the equation
\begin{equation}
\begin{split}& (\delta r)_{\lambda_{1},\lambda_{2}, \lambda_{3}} + k(\lambda_{1},\lambda_{2}, \lambda_{3})- \frac{1}{6}E(\lambda_{1},\lambda_{2}, \lambda_{3})= \\
& r_{\lambda_{2},\lambda_{3}}
 -r_{\lambda_{1}+\lambda_{2},\lambda_{3}}
 +r_{\lambda_{1},\lambda_{2}+\lambda_{3}}
 -r_{\lambda_{1},\lambda_{2}}+k(\lambda_{1},\lambda_{2},\lambda_{3}) - \frac{1}{6}E(\lambda_{1},\lambda_{2}, \lambda_{3}) = 0.
\end{split}
\end{equation}
This can be solved by
\[r_{\lambda_{1},\lambda_{2}} = 2E(\lambda_{1},i\lambda_{2},i\lambda_{1})+E(\lambda_{2},i\lambda_{2},i\lambda_{1})
\]
as is easily checked by equation (\ref{eqn:helper}).
In conclusion
\begin{equation}\label{eqn:FinalBetaPrime}
\begin{split}
& \beta'_{\lambda_{1},\lambda_{2}} =u_{\lambda_{1},\lambda_{2}} +r_{\lambda_{1},\lambda_{2}} \\&= \sum_{i<j<k}\Big{[}\frac{5}{12}
E(n_{1,i}\lambda^{i},n_{1,j}\lambda^{j},n_{2,k}\lambda^{k})
+\frac{9}{12}
E(n_{1,i}\lambda^{i},n_{2,j}\lambda^{j},n_{1,k}\lambda^{k})
+\frac{1}{12}
E(n_{2,i}\lambda^{i},n_{1,j}\lambda^{j},n_{1,k}\lambda^{k}) \\
&  \quad +
\frac{5}{12}E(n_{1,i}\lambda^{i},n_{2,j}\lambda^{j},n_{2,k}\lambda^{k})
+\frac{9}{12}E(n_{2,i}\lambda^{i},n_{1,j}\lambda^{j},n_{2,k}\lambda^{k})
+\frac{1}{12}E(n_{2,i}\lambda^{i},n_{2,j}\lambda^{j},n_{1,k}\lambda^{k})
\Big{]} \\
& \quad + 2E(\lambda_{1},i\lambda_{2},i\lambda_{1})+E(\lambda_{2},i\lambda_{2},i\lambda_{1}).
\end{split}
\end{equation}
Since $\beta'$ depends additively on $E$, we have provided the promised splitting of (\ref{eqn:GerbeSemiCharSplit}) sending $E$ to $\beta'$.

\section{Appendix 2}\label{appendix2}
In this appendix we need to explain a natural way to go from group cocycles to geometric objects: line bundles and gerbes on stacks.  When two cocycles differ by the boundary of a cochain, then this cochain can be assigned in a natural way an isomorphism between the geometric objects determined by the two cocycles.  The material in this appendix is more or less just a quick summary of certain tools from Polishchuk's recent preprint \cite{Po2008}.  We have added only minor details and changed some notation and conventions to match the philosophy of the current work.  Let a discrete group $\Gamma$ act on the right on a space $W$ so we can form the stack $[W/\Gamma]$.  Then there
is a functor
\begin{equation}\label{eqn:1CyclesToGroupoid}
[Z^{1}(\Gamma, \mathcal{O}^{\times}(W))/C^{0}(\Gamma, \mathcal{O}^{\times}(W))] \to \mathcal{P}ic([W/\Gamma])
\end{equation}
inducing the standard isomorphism of groups
\[H^{1}(\Gamma, \mathcal{O}^{\times}(W)) \cong \text{ker}[\text{Pic}([W/\Gamma])\to \text{Pic}(W)].
\]
On the level of objects, the functor (\ref{eqn:1CyclesToGroupoid}) is defined as the map
\begin{equation}\label{eqn:Ob1CyclesToGroupoid}
Z^{1}(\Gamma, \mathcal{O}^{\times}(W)) \to \textup{ob}(\mathcal{P}ic([W/\Gamma])
\end{equation}
denoted as
\[\phi \mapsto L_{\phi}
.\] Suppose we are given $\phi \in Z^{1}(\Gamma, \mathcal{O}^{\times}(W))$ and a diagram
\begin{equation}\label{eqn:StackMap}
\xymatrix{\ar @{} [dr] |{}
&P \ar[d]^-{\varpi} \ar[r]^-{c} & W &  \\
  &
 U  &   & }
\end{equation}
where $\varpi$ is a $\Gamma$-principal bundle and $c$ is $\Gamma$-equivariant.
Then let
\[L_{\phi}(U) = \{ f \in \mathcal{O}^\times(P) | \gamma \cdot f = (\phi(\gamma) \circ c)f \}.
\]
This defines an $\mathcal{O}^{\times}$-torsor on $[W/\Gamma]$.   On the level of morphisms , the functor (\ref{eqn:1CyclesToGroupoid}) is described as the map
\begin{equation}\label{eqn:Mor1CyclesToGroupoid}
I: \delta^{-1}(\phi_{2} - \phi_{1}) \to \text{Hom}(L_{\phi_{1}},L_{\phi_{2}})
\end{equation}
which we now describe.
Given a boundary $\kappa \in C^{0}(\Gamma, \mathcal{O}^{\times}(W))$ such that
\[\delta \kappa= \phi_{2}-\phi_{1}\]
we can get an isomorphism \[I_{\kappa}: L_{\phi_1} \to L_{\phi_2}\]
given by
\[f \mapsto (c^{*}\kappa) f.
\]
\begin{rmk}\label{rmk:GeomSheaves}
To get a stack over $[W/\Gamma]$ on which $\mathbb{C}^{\times}$ acts with quotient $[W/\Gamma]$, one simply considers the quotient
\[[(W \times \mathbb{C}^{\times}) / \Gamma] \to [W/\Gamma]
\]
determined by $\phi$.  Suffice it to say that one can get a functor from $[Z^{1}(\Gamma, \mathcal{O}^{\times}(W))/C^{0}(\Gamma, \mathcal{O}^{\times}(W))]$ into such stacks and following this with the functor of taking sections gives a functor naturally equivalent to the one described.
\end{rmk}

In \cite{Po2008}, Polishchuk defines a category (see Remark \ref{rmk:TwoCat}) of {\it 1-cocycles} of $\Gamma$ with values in the groupoid $\mathcal{P}ic(W)$.  He uses these cocycles to define certain interesting categories of sheaves which serve as non-commutative objects in algebraic geometry.  In fact, he actually defines a category of 1-cocycles to a group acting on any groupoid.  In particular one can use this as a tool to describe gerbes on families of complex tori.  This can be used to unify previous work of the author \cite{Be2006} with the current project.  For the current case we will only use the action groupoid of $\Gamma$ acting on $W$.  We will denote the category by $\mathcal{Z}^{1}(\Gamma, \mathcal{P}ic(W))$.  An object in this category will be denoted $(\mathcal{L}, N)$.  It consists of a collection of  $\mathcal{O}^{\times}$-torsor $\mathcal{L}_{\gamma}$ for
every $\gamma \in \Gamma$ together with isomorphisms
\begin{equation} \label{eqn:LineBundleIsoms}
N_{\gamma_{1},\gamma_{2}}: \mathcal{L}_{\gamma_{1}} \otimes \gamma_{1*}\mathcal{L}_{\gamma_{2}} \to \mathcal{L}_{\gamma_{1} \gamma_{2}}
\end{equation}
satisfying the natural consistency condition
\[N_{\gamma_{1}\gamma_{2}, \gamma_{3}}\circ (N_{\gamma_{1}, \gamma_{2}} \otimes 1)
=N_{\gamma_{1},\gamma_{2}\gamma_{3}} \circ (1\otimes \gamma_{1*} N_{\gamma_{2},\gamma_{3}})
\]
of maps
\[\mathcal{L}_{\gamma_{1}} \otimes \gamma_{1*} \mathcal{L}_{\gamma_{2}} \otimes \gamma_{1*} \gamma_{2*} \mathcal{L}_{3} \to \mathcal{L}_{\gamma_{1}\gamma_{2}\gamma_{3}}.
\]
  A morphism $(\mathcal{L}^{(1)}, N^{(1)}) \to (\mathcal{L}^{(2)}, N^{(2)})$ consists of a pair $(\mathcal{M},K)$, comprising an $\mathcal{O}^{\times}$ torsor $\mathcal{M}$ on $W$, together with isomorphisms
\begin{equation}\label{eqn:Kdef}
K_{\gamma}:\mathcal{L}^{(2)}_{\gamma} \to  \mathcal{M} \otimes \mathcal{L}^{(2)}_{\gamma} \otimes \gamma_{*} \mathcal{M}^{-1}
\end{equation}
satisfying
\[N^{(1)}_{\gamma_{1}, \gamma_{2}} \circ (K_{\gamma_{1}} \otimes \gamma_{1*} K_{\gamma_{2}}) = K_{\gamma_{1}\gamma_{2}} \circ N^{(2)}_{\gamma_{1}, \gamma_{2}}
\]
as maps
\[\mathcal{L}^{(2)}_{\gamma_{1}} \otimes \gamma_{1*} \mathcal{L}^{(2)}_{\gamma_{2}} \to \mathcal{M} \otimes \mathcal{L}^{(1)}_{\gamma_{1} \gamma_{2}} \otimes (\gamma_{1} \gamma_{2})_{*} \mathcal{M}^{-1}
.\]
We now describe a functor (see remark \ref{rmk:TwoCat}) from 1-cocycles of $\Gamma$ in $\mathcal{P}ic(W)$ to gerbes on $[W / \Gamma]$ whose pullback to $W$ is trivial:

\begin{equation}\label{eqn:Groupoid1CyclesToGerbes}
\mathcal{Z}^{1}(\Gamma,\mathcal{P}ic(W)) \to \mathfrak{G}erbes([W/\Gamma])
\end{equation}
On the level of objects
the functor (\ref{eqn:Groupoid1CyclesToGerbes}) consists of a map
\begin{equation}\label{eqn:ObGroupoid1CyclesToGerbes}\textup{ob}(\mathcal{Z}^{1}(\Gamma,\mathcal{P}ic(W))) \to \textup{ob}(\mathfrak{G}erbes([W/\Gamma]))
\end{equation}
to be denoted
\[(\mathcal{L},N) \mapsto \mathfrak{G}_{\mathcal{L},N}.
\]
Given a map $U \to [W/\Gamma]$ given by a diagram (\ref{eqn:StackMap}), we assign the groupoid $\mathfrak{G}_{\mathcal{L},N}(U)$ whose objects consist of a pair of a $\mathcal{O}^{\times}$-torsor $\Xi$ on $P$ together with a collection of isomorphisms
\[G_{\gamma}:   c^{*} \mathcal{L}_{\gamma} \otimes  \gamma_{*} \Xi \to \Xi
\]
satisfying
\[c^{*} N_{\gamma_{1}, \gamma_{2}} \otimes G_{\gamma_{1} \gamma_{2}} = G_{\gamma_{1}} \circ (1\otimes \gamma_{1*}G_{\gamma_{2}})
\]
together with the obvious notion of morphism.  We also have the map
\begin{equation}\label{eqn:MorGroupoid1CyclesToGerbes}\textup{mor}(\mathcal{Z}^{1}(\Gamma,\mathcal{P}ic(W))) \to \textup{mor}(\mathfrak{G}erbes([W/\Gamma]))
\end{equation}
\[\text{Hom}((\mathcal{L}^{(1)}, N^{(1)}) ,(\mathcal{L}^{(2)}, N^{(2)}) ) \to \text{Hom}(\mathfrak{G}_{\mathcal{L}^{(1)}, N^{(1)}} ,\mathfrak{G}_{\mathcal{L}^{(2)}, N^{(2)}} )
\]
\[(\mathcal{M},K) \mapsto I_{\mathcal{M},K}
\]
Now say we are given a morphism
\[(\mathcal{M}, K):(\mathcal{L}^{(1)}, N^{(1)}) \to (\mathcal{L}^{(2)}, N^{(2)})
\]
as above.  We need to assign this to an isomorphism of gerbes
\[\mathfrak{G}_{\mathcal{L}^{(1)}, N^{(1)}} \to \mathfrak{G}_{\mathcal{L}^{(2)}, N^{(2)}}
\]
which we provide by giving a compatable collection of groupoid isomorphisms
\[\mathfrak{G}_{\mathcal{L}^{(1)}, N^{(1)}}(U) \to \mathfrak{G}_{\mathcal{L}^{(2)}, N^{(2)}}(U)
.\]
for each $U$ as in (\ref{eqn:StackMap}).
Such an isomorphism is given by
\[\Xi \mapsto  c^{*} \mathcal{M} \otimes \Xi
\]
and
\[G^{(1)}_{\gamma} \circ K_{\gamma} = G^{(2)}_{\gamma}
\]
where the definition of the functor on the level of morphisms on the two groupoids is clear.

Let $B\mathbb{C}^{\times}$ denote the groupoid of $\mathbb{C}^{\times}$-torsors.  Given such a 1-cocycle, allows one to define the action of the group $\Gamma$ on $W \times B\mathbb{C}^{\times}$ in a way that commutes with the action of $\Gamma$ on the first factor only.  Indeed for each $\gamma \in \Gamma$ one has a functor $F_{\gamma}$ in $Aut(W \times B\mathbb{C}^{\times})$
given on objects by
\[(w, M) \to (w \gamma, \mathcal{L}_{\gamma}|_{w\gamma^{-1}} \otimes M).
\]
The definition on morphisms is clear.  Then there exist natural isomorphisms
\[F_{\gamma_{1}} \circ F_{\gamma_{2}} \Rightarrow F_{\gamma_{1} \gamma_{2}}.
\]
corresponding to (\ref{eqn:LineBundleIsoms}) and these obey a coherence condition corresponding to that mentioned above.  The result is that one
can take the quotient in a way that one still has a map to the original stack
\begin{equation} \label{eqn:GeomGerbe}
[(W \times B \mathbb{C}^{\times}) / \Gamma] \to [W / \Gamma].
\end{equation}
The sections of this projection form a $\mathcal{O}^{\times}$ gerbe isomorphic to $\mathfrak{G}_{\mathcal{L}, N}$ and in general the situation is analogous to that in Remark \ref{rmk:GeomSheaves}.

There is a functor (see remark \ref{rmk:TwoCat})
\begin{equation}\label{eqn:2CyclesToGroupoid1Cycles}
[Z^{2}(\Gamma, \mathcal{O}^{\times}(W))/C^{1}(\Gamma, \mathcal{O}^{\times}(W))] \to \mathcal{Z}^{1}(\Gamma,\mathcal{P}ic(W))
\end{equation}
which we now describe.  On objects it is the map
\begin{equation}\label{eqn:ob2CyclesToGroupoid1Cycles}
Z^{2}(\Gamma, \mathcal{O}^{\times}(W)) \to \text{ob}(\mathcal{Z}^{1}(\Gamma, \mathcal{P}ic(W))
\end{equation}
given by
\[\Phi \mapsto (\mathcal{O},\Phi)
\]
defined for $\Phi \in Z^{2}(\Gamma, \mathcal{O}^{\times}(W))$ by choosing the line bundle to be trivial and the isomorphism in (\ref{eqn:LineBundleIsoms}) to be multiplication by $\Phi_{\gamma_{1}, \gamma_{2}}$.
On the level of morphisms the map is
\begin{equation}\label{eqn:mor2CyclesToGroupoid1Cycles}
\delta^{-1}(\Phi_{2} - \Phi_{1}) \to \text{Hom}((\mathcal{O},\Phi_{1}),(\mathcal{O},\Phi_{2}))
\end{equation}
\[C \mapsto (\mathcal{O},C)
\]
here $C \in C^{1}(\Gamma, \mathcal{O}^{\times}(W))$ satisfying
\begin{equation}\label{eqn:some_name}
\Phi_{2} - \Phi_{1} = \delta C
\end{equation}
 and $(\mathcal{O}, C)$ is comprised of the isomorphisms (\ref{eqn:Kdef})  consisting of multiplication by $C_{\gamma}$.

Combining the functors (\ref{eqn:2CyclesToGroupoid1Cycles}) and (\ref{eqn:Groupoid1CyclesToGerbes}) we get a functor (see remark \ref{rmk:TwoCat})
\begin{equation} \label{eqn:2CyclesToGerbes}[Z^{2}(\Gamma, \mathcal{O}^{\times}(W))/C^{1}(\Gamma, \mathcal{O}^{\times}(W))] \to \mathfrak{G}erbes([W/\Gamma]).
\end{equation}
On objects this is the map
\begin{equation}\label{eqn:ob2CyclesToGerbes}
Z^{2}(\Gamma, \mathcal{O}^{\times}(W)) \to \text{ob}(\mathfrak{G}erbes([W/\Gamma]))
\end{equation}
to be denoted by
\[\Phi \mapsto \mathfrak{G}_{\Phi}
\]
and on morphisms
\begin{equation}\label{eqn:mor2CyclesToGerbes}
\delta^{-1}(\Phi_{2} - \Phi_{1}) \to \text{Hom}(\mathfrak{G}_{\Phi_{1}}, \mathfrak{G}_{\Phi_{2}})
\end{equation}
\[C \mapsto I_{C}
\]
\begin{rmk} \label{rmk:TwoCat}
Notice that in all three categories we have
mentioned in this appendix, \[
[Z^{2}(\Gamma, \mathcal{O}^{\times}(W))/C^{1}(\Gamma, \mathcal{O}^{\times}(W))], \ \
\mathcal{Z}^{1}(\Gamma,\mathcal{P}ic(W)), \ \  \text{and} \ \
\mathfrak{G}erbes([W/\Gamma])\] it is unreasonable to demand that
morphisms satisfy the associativity condition on the nose.  Rather,
for every three morphisms, there is a canonical natural isomorphism
between the two ways of composing the three morphisms.  This natural
isomorphism satisfies a consistency condition involving four
morphisms.  Also the functors we described are not really functors on
the nose, but rather there is a natural isomorphism (consistent with
the associativity ones) between the functor applied to a composition
of two morphisms, and the composition of the images of the two
morphisms.  There are two ways of dealing with these issues.  One way
is to force each of the three categories to become an honest category by
replacing morphisms with the obvious notion of equivalence classes of
morphisms.  With that definition we have associativity and all our functors as we have described them correspond to honest functors.  The other way is to treat each of the three categories as
a $2-$category.  While we leave most of the details to the meticulous
reader we simply remark that all the functors we have described can be
promoted to functors of $2-$ categories.  A $2-$morphism in
$[Z^{2}(\Gamma, \mathcal{O}^{\times}(W))/
C^{1}(\Gamma, \mathcal{O}^{\times}(W))]$ between
two elements $C^{(1)}$ and $C^{(2)}$ which both satisfy
$(\ref{eqn:some_name})$ is an element
$D \in C^{0}(\Gamma, \mathcal{O}^{\times}(W))$ such that
\[\delta D + C^{(1)} = C^{(2)}.
\]
A $2-$morphism in $\mathcal{Z}^{1}(\Gamma,\mathcal{P}ic(W))$ between
two pairs $(\mathcal{M}_{1},K^{(1)})$ and $(\mathcal{M}_{2},K^{(2)})$ of isomorphisms from $(\mathcal{L}^{(1)}, N^{1})$ to $(\mathcal{L}^{(2)}, N^{2})$ as in (\ref{eqn:Kdef}) consists of an isomorphism
$y: \mathcal{M}_{2} \to  \mathcal{M}_{1}$ satisfying
\[(y^{-1} \circ 1 \circ \gamma_{*} y^{\vee}) \circ K^{(1)}_{\gamma} = K^{(2)}_{\gamma}
\]
as maps
\[\mathcal{L}^{(2)}_{\gamma} \to \mathcal{M}_{2} \otimes \mathcal{L}^{(1)}_{\gamma} \otimes \mathcal{M}_{2}^{-1}.
\]
The  $2-$morphisms in $\mathfrak{G}erbes([W/\Gamma])$ are all
given by multiplication by some nowhere zero $\Gamma$-invariant
holomorphic function on $P$.  All these notions can be made to
correspond with the previously given maps on $1-$morphisms in
an obvious way.
\end{rmk}

\vskip 0.2in \noindent {\scriptsize {\bf Oren Ben-Bassat,}
Department of Mathematics, University of Haifa, Mount Carmel, 31905, ISRAEL, oren.benbassat@gmail.com}



\begin{thebibliography}{99}

\bibitem{Ap1891}
P. Appell, {\em Sur les fontions p\'eriodiques de deux variables,}
J. de Math S\'er. IV, 7, 157-219 (1891)

\bibitem{BeXu2003}
K. Behrend and P. Xu, {\em {$S\sp 1$}-bundles and gerbes over differentiable stacks,}
C. R. Math. Acad. Sci. Paris, 336, 2, 163-168 (2003), {\tt http://arxiv.org/abs/math/0306182}

\bibitem{Be2006}
O. Ben-Bassat, {\em Twisting derived equivalences,}
Trans. Amer. Math. Soc.,  361, 5469-5504 (2009), {\tt http://arxiv.org/abs/math/0606631}

\bibitem{BeBlPa2007}
O. Ben-Bassat, J. Block, and T. Pantev, {\em Non-commutative tori and {F}ourier-{M}ukai duality,} Compositio Mathematica, 143, 423-475 (2007), {\tt http://arxiv.org/abs/math/0509161}

\bibitem{Be1972}
V. G. Berkovi{\v{c}}, {\em The {B}rauer group of abelian varieties,}
Funkcional Anal. i Prilo\v zen. (6), 3, 10-15 (1972)

\bibitem{BiLa1999}
C. Birkenhake and H. Lange, {\em Complex Tori,} Progress in mathematics v. 177, Birkhauser (1999)

\bibitem{BiLa2004}
C. Birkenhake and H. Lange, {\em Complex abelian varieties,} Fundamental
Principles of Mathematical Sciences, 302, Second Edition, Springer-Verlag (2004)


\bibitem{Br1993}
J.-L. Brylinski, {\em Loop spaces, characteristic classes
and geometric quantization,} Progress in Mathematics,
107, Birkh\"auser Boston Inc. (1993)

\bibitem{De}
J.-P. Demailly, {\em Complex analytic and algebraic geometry,} OpenContent Book {\tt http://www-fourier.ujf-grenoble.fr/\~{}demailly/books.html}

\bibitem{Di2004}
A. Dimca, {\em Sheaves in Topology,} Universitext, Springer-Verlag (2004)

\bibitem{DoGa2002}
R. Dongai and D. Gaitsgory, {\em The gerbe of Higgs bundles,} Transform. Groups,
7, 109-153 (2002), {\tt http://arxiv.org/abs/math/0005132}

\bibitem{DoPa2008}
R. Donagi and T. Pantev, {\em Torus fibrations, {G}erbes, and duality, With an appendix by Dmitry Arinkin,} Mem. Amer. Math. Soc. (193), 901 (2008), {\tt http://arxiv.org/abs/math/0306213}

\bibitem{ElNa1983}
G. Elencwajg and M. S. Narasimhan, {\em Projective bundles on a complex torus,}
J. Reine Angew. Math. (340), 1-5 (1983)

\bibitem{FeHeRoZh2008}
G. Felder, A. Henriques, C. Rossi, and C. Zhu, {\em A gerbe for the elliptic gamma function,} Duke Math. J. (141), 1, 1-74 (2008), {\tt http://arxiv.org/abs/math/0601337}


\bibitem{Gi1971}
J. Giraud, {\em Cohomologie non ab\'elienne,} Die Grundlehren der mathematischen Wissenschaften, Band 179 (1971)


\bibitem{Ho1972}
R. Hoobler, {\em Brauer groups of abelian schemes,}
Ann. Sci. \'Ecole Norm. Sup. (4), 5, 45-70 (1972)

\bibitem{Hu1893}
G. Humbert, {\em Theorie generale des surfaces hyperelliptiques,}
J. de Math S\'er. IV, 9, 29-170 and 361-475 (1893)

\bibitem{Mu1970}
D. Mumford, {\em Abelian varieties,} Tata Institute of Fundamental Research Studies in Mathematics, No. 5 (1970)

\bibitem{Po2003}
A. Polishchuk, {\em Abelian varieties, theta functions and the {F}ourier
transform,} Cambridge Tracts in Mathematics, 153, Cambridge University Press (2003)

\bibitem{Po2008}
A. Polishchuk, {\em Kernel algebras and generalized Fourier-Mukai transforms,} {\tt http://arxiv.org/abs/0810.1542}

\bibitem{Sc2005}
S. Schroeer, {\em Topological methods for complex-analytic Brauer groups,} Topology, 44, no. 5, 875-894 (2005),
{\tt http://arxiv.org/abs/math/0405223}

\bibitem{We1958}
A. Weil, {\em Introduction \`a l'\'etude des vari\'et\'es k\"ahl\'eriennes,}
Publications de l'Institut de Math\'ematique de l'Universit\'e de
Nancago, VI. Actualit\'es Sci. Ind. no. 1267 (1958)



\end{thebibliography}
\end{document}